%
%

%
%


\input amstex
\documentstyle{amsppt}
\NoBlackBoxes
\NoRunningHeads
\TagsOnLeft

\hsize = 5.5 truein
\vsize = 8.5 truein

\define\vre{\varepsilon}
\define\hd{Hausdorff dimension}
\define\hs{homogeneous space}
\define\df{\overset\text{def}\to=}
\define\un#1#2{\underset\text{#1}\to#2}
\define\br{\Bbb R}
\define\bn{\Bbb N}
\define\bz{\Bbb Z}
\define\bq{\Bbb Q}

\define\di{Diophantine}
\define\da{Diophantine approximation}
\define\de{Diophantine exponent}
\define\va{\bold a}
\define\ve{\bold e}
\define\vh{\bold h}
\define\vx{\bold x}
\define\vy{\bold y}

\define\vv{\bold v}
\define\vu{\bold u}
\define\vr{\bold r}

\define\vw{\bold w}
\define\vp{\bold p}
\define\vq{\bold q}

\define\vc{\bold c}

\define\vf{\bold f}

\define\vg{\bold g}

\define\rk{\operatorname{rk}}

\define\spr{Sprind\v zuk}

\define\nz{\smallsetminus \{0\}}

\define\be{Besicovitch}

\define\cag{$(C,\alpha)$-good}

\define\p{\Phi}
\define\vrn{\varnothing}
\define\ssm{\smallsetminus}

\define\mr{M_{m,n}}
\define\amr{$A\in M_{m,n}$}

\define\SL{\operatorname{SL}}
\define\GL{\operatorname{GL}}
\define\Id{\operatorname{Id}}
\define\const{\operatorname{const}}

\define\supp{\operatorname{supp}}

\topmatter
\title 
An extension of quantitative nondivergence and applications to 
Diophantine exponents
\endtitle  

\author { Dmitry Kleinbock} \\ 
\endauthor

    \address{ Dmitry Kleinbock,  Department of
Mathematics, Brandeis University, Waltham, MA 02454-9110}
  \endaddress

\email kleinboc\@brandeis.edu \endemail

  \thanks 
This work was supported in part by NSF
Grant DMS-0239463.
\endthanks       


\abstract 
We present a sharpening of nondivergence
estimates for unipotent (or more generally polynomial-like) flows on \hs s.
Applied to metric \da, it yields  precise formulas for \de s of
affine subspaces of $\br^n$ and their nondegenerate submanifolds.

\endabstract



\endtopmatter
\document

\heading{0. Introduction}\endheading


This paper continues the theme started  in 1971  by G.\,A.\ Margulis \cite{Mr1} 
when he showed
that trajectories of one-parameter 
unipotent flows on $\SL_k(\br)/\SL_k(\bz)$  are never divergent.
This had been earlier conjectured by I.\ Piatetski-Shapiro, and was used by
Margulis for the proof of the Arithmeticity Theorem for nonuniform lattices.
A decade later, 
   S.\,G.\ Dani  \cite{Da1--2, Da4} modified the method of Margulis, showing
 that any unipotent orbit returns to big compact subsets with high
frequency. 
The latter statement was part
of  Dani's proof \cite{Da1} of  finiteness of  locally  finite 
unipotent-invariant ergodic measures  and was in M.\ Ratner's  proof  \cite{Rt1--2} of
Raghunathan's topological conjecture.

The next development came in 1998, when a very general explicit estimate for the
above frequency in terms of the size of compact sets was given in the paper
of Margulis and the 
author \cite{KM}. In fact it was
done in a bigger generality, namely for
a large class of maps
from $\br^d$ into $\SL_k(\br)$, which made it possible to derive important
applications  to metric \da\ on manifolds.

To state some of the results from that paper, recall that 
the space $$\Omega_k\df \SL_k(\br)/\SL_k(\bz)$$  can be identified with
the  space of  unimodular
lattices in $\br^k$, and that
$\Omega_k = \bigcup_{\vre > 0} K_\vre$, where the sets
$$K_\vre \df \big\{\Lambda \in \Omega_k \bigm| \|\vv\| \ge \vre\text{ for all
}\vv\in\Lambda\nz\big\}
$$ are compact  
(Mahler's Compactness Criterion; see e.g.\ 
\cite{Rg, Corollary
10.9}). Here
 $\|\cdot\|$ can be any norm on  $\br^k$, which we will assume to be Euclidean, and
extend to the space of  discrete 
subgroups of   $\br^k$ by letting $\|\Gamma\|$  be  
the volume of the quotient space $\Gamma_\br/\Gamma$
(here and hereafter $\Gamma_\br$ stands for the $\br$-linear span of $\Gamma$)
if $\Gamma
\ne\{0\}$, and $1$ otherwise.

\medskip

Another notion 
we need to introduce is that of  functions
\cag\ on an open subset of $\br^d$. We postpone a definition of this property
until the next section, noting that, roughly speaking, it can be interpreted as some
kind of polynomial-like behavior. 

\medskip

The following is a special case of  one of the main theorems from \cite{KM}:
 \proclaim{Theorem  0.1} Let $d,k\in\bn$, $C,\alpha > 0$, $0 < \rho  \le 1$,  and
let  a ball $B = B(\vx_0,r) \subset \br^d$ and a map $h:\tilde B \df
B\big(\vx_0,3^{k}r\big) \to \SL_k(\br)$  be given.  
 Assume that the following two conditions hold:
$$
\text{for any 
subgroup $\Gamma\subset\bz^k$, the function
$\vx\mapsto\|h(\vx)\Gamma\|$ is \cag\ on $\tilde B$,}\tag 0.1
$$
and
$$
\text{for any subgroup }
\Gamma\subset\bz^k, \quad\sup_{\vx\in
B}\|h(\vx)\Gamma\| \ge \rho\,.\tag 0.2
$$
Then 
$$
\forall\, 0 < \vre \le \rho\quad\lambda\big(\{\vx\in B\mid h(\vx)\bz^k \notin K_\vre
\}\big)
\le
\const(d,k)\cdot C\left(\frac\vre
\rho \right)^{\alpha}  \lambda(B)\,,\tag 0.3
$$
where the constant above is explicitly computable and depends only on $d$ and $k$. 
\endproclaim
 
Here and hereafter $\lambda$ stands for  Lebesgue measure on  $\br^d$.

\medskip

Roughly speaking, the informal meaning of this theorem is as follows. 
Under the presence of condition (0.1), which captures the polynomial-like behavior of the
map
$h$, one of the following two alternatives holds: either most of the `orbit'
$h(B)\bz^k$ is contained  in $K_\vre$, or there exists $\Gamma\subset\bz^k$ which is
`responsible for the whole orbit being far away', namely, such that the covolume of 
$h(\vx)\Gamma$ is uniformly small for all $\vx\in B$.

Note that in the papers of Margulis and Dani the function
$h$ was of  the form $x\mapsto u_xg$, where $\{u_x\}$ was a unipotent subgroup of
$\SL_k(\br)$ and $g$ a fixed element of $\SL_k(\br)$.
Another class of important applications of Theorem  0.1 is to \da. There, to study 
\di\ properties of almost every vector of the form $\vf(\vx)$, where $\vf$ is a map from
$\br^d$ to $\br^n$, one takes
$k = n + 1$ and considers $$
h(\vx) = \left(\matrix
e^{-t} & 0  \\
0 & e^{t/n}I_n
\endmatrix \right)\left(\matrix
1 & \vf(\vx)
\\
0 & I_n
\endmatrix \right)\,.\tag 0.4
$$ 
In this case checking (0.1) amounts to establishing polynomial-like  behavior of $\vf$,
and (0.2) can often be extracted from global  \di\ properties of the image
of
$\vf$. The papers \cite{KM, BKM, BBKM, G1, K2--3, KLW, KT, KW1--3} contain various number-theoretic
applications of Theorem  0.1 and its generalizations. See also \cite{K1, Mr2}
for a   survey of this method and related results.

The purpose of the present paper is to pay special attention to condition (0.2) of
Theorem  0.1. Namely,  it is not surprising that
the higher is the {\sl rank\/}  $\rk(\Gamma)$ of
$\Gamma$ (defined as the rank of $\Gamma$ as a $\bz$-module, or,
equivalently, the dimension of 
$\Gamma_\br$),
 the harder is usually the task
of estimating
$\|h(\cdot)\Gamma\|$ from below. For example if   $\Gamma$ is of rank $j$ and
$h(\vx)\Gamma$ is 
generated by  $j$ orthogonal vectors of length $\rho$, one has $\|h(\vx)\Gamma\| \asymp
\rho^{j}$. 
Indeed, in some cases relevant
to
\di\ applications one can prove estimates of the following form:
$$
\text{for any }
\Gamma\subset\bz^k,\quad\sup_{\vx\in
B}\|h(\vx)\Gamma\| \ge \rho^{\rk(\Gamma)}\,.\tag 0.2$'$
$$
This leads to a natural question -- whether
or not such a condition is enough for deriving (0.3). It became clear to the author in
the process of working on the  paper \cite{K2} that an
affirmative answer to the latter question would make it  possible to significantly 
generalize several key results from that paper.

\medskip

The following theorem provides such an answer:

\proclaim{Theorem  0.2} Let $d,k\in\bn$, $C,\alpha,\rho > 0$, 
and
suppose that a ball $B = B(\vx_0,r) \subset \br^d$ and a map $h:\tilde B \df
B\big(\vx_0,3^{k}r\big) \to \SL_k(\br)$  satisfy {\rm (0.1)} and  {\rm (0.2$'$)}.  
Then  {\rm (0.3)}  holds. 
\endproclaim
 
We will comment on the geometric significance of the improvement of Theorem 0.2
over Theorem 0.1 at the end of \S 2, after a more general result is proved.

\medskip

Let us now describe some number-theoretic applications of the above  theorem.
For  $m,n\in\bn$,  $\mr$ will stand for the space of real matrices with  $m$ rows and $n$
columns. 
Define the {\sl
\de\/} $\omega(A)$ of \amr\ (sometimes called `the exact order' of $A$) 
 to be the supremum of $v > 0$ for which there are
infinitely many $\vq\in
\bz^n$ such that
$$
 \|A\vq + \vp\|   < \|\vq\|^{-v}  \tag 0.5
$$
for some
$\vp\in\bz^m$.  Note that this quantity is independent of the choice of
norms on $\br^m$ and $\br^n$ (thus we will repeatedly switch between
 Euclidean and supremum norms whenever it is convenient).
%
%
 It is well known and easy to see that one has $\,n/m \le \omega(A) \le
\infty$ for all
\amr, 
and $
\omega(A) = n/m$ for
$\lambda$-almost every
\amr.

We will specialize to the case $m = 1$, that is, consider \di\ properties
of  $\vy\in\br^n$ interpreted as row vectors (see however \S 6.3 for the column vector set-up). Further, our emphasis will be on \da\ 
with dependent quantities, where the dependence is expressed by means of a Borel measure
$\mu$ on $\br^n$. 
Namely, 
let us define the {\sl
\de\/}
$\omega(\mu)$ of $\mu$ to be the
$\mu$-essential supremum of the function $\omega(\cdot)$; in other words, 
$$
\omega(\mu)\df \sup \big\{\,v\bigm| \mu(\{\vy\mid \omega(\vy) > v \}) > 0\big\} 
\,.
$$
Clearly it only depends on the measure class of $\mu$.
If $\mu$ is
naturally associated with a subset $\Cal M$ of $\br^n$ supporting $\mu$ (for example,
if $\Cal M$ is a smooth submanifold of $\br^n$ and $\mu$ is 
 the measure class of the Riemannian volume on $\Cal M$ $\iff$ the pushforward
$\vf_*\lambda$ of
$\lambda$ by any smooth map $\vf$ parametrizing
$\Cal M$), we will define the \de\ $\omega(\Cal M)$ of $\Cal M$ to be equal
to that of
$\mu$. Clearly
$\omega(\mu)
\ge n$ for any
$\mu$, and $ \omega(\lambda)= \omega(\br^n)$ is equal to $n$. 
The latter justifies the terminology introduced to \da\ on manifolds
by V.\ Sprind\v zuk:  a measure $\mu$ on $\br^n$ (resp., a submanifold
$\Cal M$ of $\br^n$) is  {\sl
extremal\/} if
$\omega(\mu)$ (resp., 
$\omega(\Cal M)$) is equal to $n$, that is, 
attains the smallest
possible value.

\comment
For a map $\vf:U\to \br^n$, $U\subset
\br^d$, it is natural to define the {\sl \de\/} 
$\omega(\vf)$ of $\vf$ to be the $\lambda$-essential  infimum of
$\omega\big(\vf(\cdot)\big)$, i.e.
$$
\omega(\vf)\df \sup\left\{v\Bigm| \lambda\big(\{\vx\in U\mid \vf(\vx)\text{ is 
$v$-approximable}\}\big) >
0\right\}\,.
$$
If $\Cal M$ is a smooth manifold, we let the \de\ 
$\omega(\Cal M)$ of $\Cal M$ to be the \de\ of its
parametrizing map.
A manifold/map is said to be {\sl extremal\/} iff its \de\ is the smallest
possible, i.e.~is equal to $n = \omega(\br^n)$. 
\endcomment

It was conjectured by  \spr\ 
in 1980 \cite{Sp} and proved in \cite{KM} that real analytic manifolds not
contained in any proper affine subspace of $\br^n$ are extremal. 
More generally, let us say that a differentiable
map $\vf:U\to \br^n$, where $U$ is an open subset
of $\br^d$,  is {\sl nondegenerate in  an affine subspace $\Cal L$ of $\br^n$ at
\/}
$\vx\in U$ if $\vf(U)\subset \Cal L$ and 
the span of all the partial
derivatives of
$\vf$ at
$\vx$ up to some order coincides with the linear part of  $\Cal L$. If $\Cal M$ is
a $d$-dimensional 
submanifold of $\Cal L$, we will say
that $\Cal M$ is {\sl nondegenerate  in $\Cal L$ at $\vy\in \Cal M$} if any
(equivalently, some)  diffeomorphism
$\vf$ between an open subset $U$  of $\br^d$ and a neighborhood of
$\vy$ in $\Cal M$ is   nondegenerate  in $\Cal L$ at $\vf^{-1}(\vy)$. We will say
that
$\vf:U\to \Cal L$ (resp., $\Cal M\subset \Cal L$) is 
{\sl nondegenerate  in
$\Cal L$\/} if it is nondegenerate  in
$\Cal L$ at
$\lambda$-a.e.\ point of
$U$  (resp., of $\Cal M$, in the sense of the 
smooth measure class on $\Cal M$). 

\comment
{\sl $\ell$-nondegenerate in  an affine subspace $\Cal L$ of $\br^n$ at
\/}
$\vx\in U$ if $\vf(U)\subset \Cal L$ and 
$$
\aligned
\text{the linear part of  $\Cal L$ is spanned by}\\ \text{ partial
derivatives of
$\vf$ at
$\vx$ of order up}&\text{ to }\ell\,.
\endaligned\tag 0.4
$$  
If $\Cal M$ is
a $d$-dimensional 
submanifold of $\Cal L$, we will say
that $\Cal M$ is {\sl $\ell$-nondegenerate  in $\Cal L$ at $\vy\in \Cal M$} if any
(equivalently, some)  diffeomorphism
$\vf$ between an open subset $U$  of $\br^d$ and a neighborhood of
$\vy$ in $\Cal M$ is   $\ell$-nondegenerate  in $\Cal L$ at $\vf^{-1}(\vy)$. We will say
that
$\vf:U\to \Cal L$ (resp., $\Cal M\subset \Cal L$) is {\sl nondegenerate  in
$\Cal L$\/} at a point if (0.4) holds for some $\ell$, and  that it is 
{\sl nondegenerate  in
$\Cal L$\/} if it is nondegenerate  in
$\Cal L$ at
$\lambda$-a.e.\ point of
$U$  (resp., of $\Cal M$, in the sense of the 
smooth measure class on $\Cal M$).
\endcomment


\medskip

The next theorem generalizes some of the results of \cite{KM} and \cite{K2}:

 \proclaim{Theorem  0.3} Let $\Cal L$ be an  affine subspace  of $\br^n$,
 and let $\Cal M$ be a submanifold  of $\Cal L$ which is
nondegenerate  in
$\Cal L$. Then $$\omega(\Cal M) =
\omega(\Cal L) = 
\inf\{\omega(\vy)\mid \vy\in\Cal L\} = \inf\{\omega(\vy)\mid \vy\in\Cal M\}\,.$$
\comment
$v\ge n$. Then the following conditions are
equivalent:

\roster
\item"{\rm (i)}" $\omega(\Cal L) \le v$;
\item"{\rm (ii)}" $\omega(\Cal M) \le v$ for any submanifold $\Cal M$ of $\Cal L$
nondegenerate  in
$\Cal L$;
\item"{\rm (iii)}" $\omega(\vy) \le v$ for some $\vy\in\Cal L$.
\endroster 
\endcomment
\endproclaim

In short, \de s of affine subspaces are inherited by their nondegenerate submanifolds.
This was proved in  \cite{KM} for 
$\Cal L = \br^n$, and  in  \cite{K2} for
 extremal  $\Cal L$ (that is, with  $\omega(\Cal L)= n$). 
\comment
It is obvious
that (ii)$\Rightarrow$(i)$\Rightarrow$(iii), but the implication (iii)$\Rightarrow$(ii)
is in our opinion quite remarkable; even the weaker one, (iii)$\Rightarrow$(i),
does not seem to be known previously.

Note
that it follows from Theorem  0.3 that 
$\omega(\Cal L) =
\omega(\Cal M)$ for any
$\Cal M$ nondegenerate  in
$\Cal L$, and $$\omega(\Cal M) = \inf\{v\mid \exists\,\vy\in \Cal M\text{ which is not
$v$-approximable}\}$$ whenever 
$\Cal M$ is nondegenerate  in some affine subspace. 
\endcomment
Note that the middle equality is trivially satisfied for $\Cal L = \br^n$,
but is not at all obvious for proper
subspaces $\Cal L$. Indeed, it states that the infimum of $\omega|_\Cal L$ coincides with its 
essential supremum; that is, the existence of a single point $\vy\in\Cal
L$ with $\omega(\vy) \le v$ forces the set $\{\vy\in\Cal L\mid\omega(\vy) \le v\}$ to
have
full measure.  

Another natural problem addressed in the paper is computing
 \de s of  affine subspaces in terms of the coefficients of their parametrizing maps.
If  $\Cal L$ is an  
${s}$-dimensional affine subspace  of $\br^n$, 
by permuting variables one can 
 without loss of generality 
choose 
a parametrizing map  of the form
$\vx\mapsto(\vx,
\vx A' + \va_0)$,
where 
$A'\in M_{s,n-s}$ 
and $\va_0\in\br^{n-{s}}$ (here both $\vx$ and $\va_0$ are
row vectors). It will be convenient to denote    the matrix
$\left(\matrix
\va_0 \\ A'
\endmatrix\right)$ by $A\in M_{{s}+1,{n-{s}}}$, so that $\Cal L$ is  parametrized by $$
\vx\mapsto(\vx,
\tilde\vx A)\,,
\tag 0.6
$$
where $\tilde
\vx$ stands for $ (1,\vx)$. 

One of the advantages of such a parametrization is a possibility to relate 
\di\ properties
of $A$ to those of points of $\Cal L$. Indeed, following \cite{K2}, 
it can be easily shown that a good approximation to $A$ gives rise to a good approximation
to all points of $\Cal L$ simultaneously (see Lemma 5.4).
Consequently, $\omega(A)$ is a lower bound
for $\omega(\Cal L)$; thus
$$
\omega(\Cal L)\ge \max\big(\,\omega(A),n\big)\,.\tag 0.7$\ge$
$$
Estimating $
\omega(\Cal L)$ from above is a more difficult task. We accomplish it in
 \S 5 by writing  a precise expression for $\omega(\Cal L)$ parametrized as in 
(0.6) in terms of $A$ (see Corollary 5.2). In particular, we prove

 \proclaim{Theorem  0.4} If a proper  affine subspace $\Cal L$   of $\br^n$ is parametrized 
as in {\rm (0.6)}, where  either {\rm (a)} all the  columns or {\rm (b)} all the  rows  of $A$ 
are rational multiples of  one column (resp., row), then  {\rm (0.7$\ge$)} turns into equality; that is, one has 
$$
\omega(\Cal L)= \max\big(\,\omega(A),n\big)\,.\tag 0.7$=$
$$
\endproclaim

This generalizes \cite{K2, Theorem 1.3}. Examples of subspaces described by Theorem  0.4 include: those parallel to
coordinate subspaces; lines passing through the origin; subspaces of codimension one.
 Whether or not  {\rm (0.7$=$)} holds in general  is an open question, see \S 6.1--2 for discussion.

\medskip

We remark  that the main results of this paper, Theorems 2.1 and 2.2,
are much  more general than Theorem  0.2. Namely, we consider  maps from  \be\ metric
spaces equipped with  Federer measures; see
\S 1 for definitions. Thus our main \di\ result, Theorem 1.3, is substantially 
 more general than
Theorem  0.3. In particular, its framework includes fractal subsets of
$\br^n$ 
 or, more generally,  measures of the form
$\vf_*\mu$ where $\mu$ satisfies certain decay conditions, as in
\cite{KLW}. 
In \S 1  we review all the necessary terminology and background facts,
and in \S 2 prove the general quantitative nondivergence estimates. Section 3 describes
a connection between \da\ and dynamics which makes it possible to apply Theorem 2.2
to \de s. Theorem 1.3 is proved in \S 4 and Theorem 0.4 in \S 5. In fact, the
  \de s of  arbitrary affine subspaces are expressed in terms of so-called {\sl higher order
\de s\/} of matrices, which  are  introduced and studied in detail in \S 5.
The 
 last section 
contains
several open questions and further generalizations of the \di\ problems considered in the paper.

\medskip

\noindent {\bf Acknowledgements:}
 The author is thankful to Gregory Margulis, Barak
Weiss and the reviewer
for their valuable comments.
 This work was supported
in part by NSF
Grant DMS-0239463.

\heading{1. Preliminaries  (\be, Federer, good, nonplanar)\\
 and the main \di\
results
}
\endheading


A metric space $X$ is called {\sl   $N$-\be\/} 
if for any bounded subset $A$
of
$X$ and  any family $\Cal B$ of nonempty open balls in $X$ such that
each $x\in A$ is a center of some ball of $\Cal B$, there is a finite or countable
subfamily $\{B_i\}$ of $\Cal B$ covering $A$ with
multiplicity 
at
most $N
$.  We will say that $X$ is {\sl   \be\/} if it is  $N$-\be\ for some $N$.
 The fact that $\br^d$  is Besicovitch is the content
of Besicovitch's Covering Theorem \cite{Mt, Theorem 2.7}. 

\medskip

Let $\mu$ be a Radon
 measure
on
$X$, and $U$ an open subset of $X$ with $\mu(U) > 0$. Following \cite{KLW}, let us say that
$\mu$ 
is {\sl $D$-Federer on 
\/} $U$
        if 
$$
\sup\Sb x\in \supp\mu,\,r > 0\\ B(x,3r)\subset U\endSb\ 
\frac{\mu\big(B(x,3r)\big)}{\mu\big(B(x,r)\big)} < D\,.
\tag 1.1
$$

Equivalently, one can replace `$3$' in (1.1) by any $c > 1$, appropriately
changing the value of $D$. This explains why Federer measures are often called {\sl
doubling\/}. 
It will be useful to have a nonuniform version of the above
definition: 
 we will say that $\mu$ as above  is {\sl Federer\/} if for
$\mu$-a.e.\ $x\in X$  there
exists a neighborhood
$U$ of
$x$ and $D > 0$  such that $\mu$ is
        $D$-Federer on $U$.
Many natural examples of measures, including those supported on fractals, 
can be shown to be Federer.

\medskip

For a
subset
$B$ of
$X$ and a
function $f$ from $X$ to a normed space with norm $\|\cdot\|$, 
we let
$\|f\|_{B} \df \sup_{x\in B}\|f(x)\|$.
If  $\mu$ is a  
measure  on $X$ and
$B$ is a subset of  $X$ with $\mu(B)>0$, 
we define
         $\| f \|_{\mu,B}$ to be equal to $\|f\|_{B\,\cap\,\supp\mu}$.
\medskip

A function $f:X\to \br$ is called {\sl \cag\ on $U\subset X$ with respect to
$\mu$\/} if for any open ball
$B\subset U$ centered in $\supp\mu$ one has
$$
\forall\,\vre > 0\quad \mu\big(\big\{ x\in B\bigm|
|f(x)| < \vre\}\big) \le C\left(\frac\vre{\|f\|_{\mu,
B}}\right)^\alpha  \mu(B)\,.
\tag 1.2
$$
Informally speaking, a function is   good if the set of points
where it takes small values has small measure. 
We refer the reader to \cite{KM, BKM, KLW,  KT} for various properties and
examples. One of the elementary observations is conveniently stated below:

\proclaim{Lemma  1.1 {\rm \cite{KLW, Lemma 4.1}}}  Suppose that   $f_1,\dots,f_\ell$
are $(C,\alpha)$-good on
$U$ with respect to
$\mu$. Then $(f_1^2 + \dots +f_\ell^2)^{1/2}$ is $(\ell^{\alpha/2}C,\alpha)$-good
on
$U$ with respect to
$\mu$.\endproclaim

In the
situations when
$U$ is a  subset of 
$
\br^d$ and 
$\mu =
\lambda$,
we will omit the reference to the measure and will simply say 
`$f$ is
\cag\
on $U$', as has already been done in (0.1).  In that case
 one can replace
$\|f\|_{\mu,
B}$ in (1.2) by $\|f\|_{B}$ and not pay attention to the
restriction of the center of $B$ lying  in the support of the measure.

\medskip

\cag\ functions often come in families. For example, condition (0.1) used in 
Theorems 0.1 and 0.2 asserts that functions of the form $\vx\mapsto\|h(\vx)\Gamma\|$,
where $\Gamma$ runs through  subgroups of $\bz^k$, are all \cag\ with uniform 
$C$ and $\alpha$. Often we will need to check the \cag\ property for functions from
a given finite-dimensional function space. An example:
polynomials in $d$ variables are  \cag\ on $\br^d$ with 
$C$ and $\alpha$ depending only on $d$ and the degree of the polynomial. 

The following definition was introduced in \cite{K2} and \cite{KT}. 
Let $\vf = (f_1,\dots,f_n)$ be 
a  
map from a metric space $X$  to $\br^n$ and $\mu$ a measure on $X$. We will say that
 $(\vf,\mu)$ is {\sl good  at  $x\in
X$\/}  if there exists a neighborhood
$V$ of $x$ and positive $C,\alpha$
such that
       any linear combination of
$1,f_1,\dots,f_n$  is
$(C,\alpha)$-good on
$V$ with respect to  $\mu$. 
We will simply say  that  $(\vf,\mu)$ is 
{\sl good\/} if it is good at
$\mu$-almost every point.
Again, the reference to the measure
will be omitted
when $\mu = \lambda$, in which case we will say that $\vf$ is good or good at $x$.
For example, we will say that  polynomial
maps are  good (in fact, good at every point).
More generally, based on the work done in \cite{KM}, the following was proved in
\cite{K2}:

\proclaim{Lemma  1.2} 
Let $\Cal L$ be an affine
subspace of $\br^n$ and let  $\vf$ be a
  smooth
map  from  open  $U\subset\br^d$ to $\Cal L$ which is nondegenerate in $\Cal L$ at $\vx\in U$.
Then
$\vf$ is good at $\vx$.
\endproclaim

For a subset $M$ of $\br^n$, define its {\sl affine span\/} $\langle M\rangle_a$
to be the intersection of all affine subspaces  of $\br^n$ containing $M$.
Then  it is easy to see that $\Cal L$ in the above lemma 
is equal to $\langle \vf(B)\rangle_a$ for some open
$B\ni
\vx$.
It will be useful to define a similar property 
for more general maps and measures. Namely,  let $X$
be a   metric space,
$\mu$ a  measure on
$X$, $\Cal L$  an affine
subspace of $\br^n$ and $\vf$  a map from $X$ into $\Cal L$. 
Say that  $(\vf,\mu)$
is {\sl  nonplanar in $\Cal L$ 
\/} (cf.\ \cite{KT, \S 4} and \cite{KW3, \S 1})
   if 
$$\Cal L = \langle\vf(B\cap 
\supp\mu)\rangle_a \quad\forall
\text{  nonempty open }B\text{ with }\mu(B) > 0
\,.\tag 1.3
$$ 
As before, we will omit the dependence on $\mu$
(resp., $\Cal L$)
 when 
$\mu = \lambda$ (resp., $\Cal L = \br^n$). Clearly 
$(\vf,\mu)$ is
 nonplanar  
iff
for any nonempty open $B$ of positive measure
the
restrictions of
$1,f_1,\dots,f_{n}$ to
$B\,\cap \,\supp\,\mu$
are linearly independent over $\br$. 

\medskip

\comment
It is obvious that the three properties of maps $\vf$ considered above (nondegeneracy,
goodness, nonplanarity) are preserved under post-compositions with affine isomorphisms. 
More precisely,  let $\vf:X\to \Cal L\subset\br^n$ be $\mu$-good (resp., $\mu$-nonplanar
in $\Cal L$, nondegenerate 
in $\Cal L$) at $x\in X$ (for the last option $X = \br^d$ and $\vf$  is smooth), 
and let $\vg: \br^n\to\br^m$ be an affine map which is one-to-one on $\Cal L$;
then $\vg\circ\vf$ is also  $\mu$-good (resp., $\mu$-nonplanar
in $\vg(\Cal L)$, nondegenerate 
in $\vg(\Cal L)$) at $x$. In particular, nondegeneracy/nonplanarity 
in $\Cal L$ can be expressed in terms of nondegeneracy/nonplanarity 
of certain post-compositions in $\br^k$, where $k = \dim(\Cal L)$. This is precisely how
Lemma 1.1 is deduced from the results of \cite{KM}.

A basic example is given by nondegenerate smooth maps from $\br^d$ to $\br^n$ and $\mu = \lambda$:
 it is
clear from the  definition that nondegeneracy in $\Cal L$ implies nonplanarity in $\Cal
L$. Note that for  real analytic maps even more can be said \cite{K2}:

\proclaim{Lemma  1.2}  Let $\vf$ be a
real analytic map from a connected open subset $U$ of
$\br^d$ to $\br^n$. Then there exists an affine subspace $\Cal L$ 
of $\br^n$ which is equal to the affine span of $\vf$ at every point of
$U$, and  $\vf$ is nondegenerate
in $\Cal L$ at every point of
$U$. Consequently, $\vf$ is good and nonplanar in $\Cal L$ at every point of
$U$.
\endproclaim

\endcomment

As was said before, a basic example is given by nondegenerate smooth maps from $\br^d$
to
$\br^n$: it is clear from the  definition that nondegeneracy in $\Cal L$ implies
nonplanarity in
$\Cal L$. Thus the following statement generalizes 
Theorem  0.3:

 \proclaim{Theorem 1.3}  Let  $\mu$ be a Federer measure on a \be\ metric space $X$,
$\Cal L$ an affine subspace of
$\br^n$, and let $\vf:X\to\Cal L$ be a continuous map such that $(\vf,\mu)$ is good
and
nonplanar in $\Cal L$.
Then $$\omega(\vf_*\mu)
= \omega(\Cal L) = 
\inf\{\omega(\vy)\mid \vy\in\Cal L\} = \inf\big\{\omega\big(\vf(x)\big)\bigm| x\in\supp\,\mu\big\}
\,.
$$
\endproclaim

The special case $\Cal L = \br^n$ was (in a slightly different terminology)
one of the main results of \cite{KLW}. Note that in all the applications considered in this paper we will take
$X$ to
be an open subset of $\br^d$; however one can also work with vector spaces
over other local fields 
and, using methods from \cite{KT} and \cite{G2}, obtain 
non-Archimedean version of many results from the present
paper; see \S 6.6 for further discussion.

\medskip

Many nontrivial examples of measures $\mu$ and maps $\vf$ satisfying 
the conditions of the above theorem can be found in the paper \cite{KLW}. For example,
it is not hard to see that a measure $\mu$ on $\br^n$ is  {\sl friendly\/} 
(a property introduced in  \cite{KLW}) iff it is
Federer and $(\Id,\mu)$ is good
and
nonplanar (here $\Id$ is the identity map $\br^n\to\br^n$). Many measures naturally arising from geometric constructions 
can be shown to possess an even stronger property -- such measures were referred to as 
`absolutely decaying and  Federer' in  \cite{KLW} and as `absolutely friendly'
 in \cite{PV}; many examples of those can be found in \cite{KLW, U, SU}. It was proved in  
\cite{KLW, \S 7} that if $\mu$ is absolutely decaying and  Federer and $\vf$ is nondegenerate at
$\mu$-a.e.\  point of $\br^d$,
then $(\vf,\mu)$ is good
and
nonplanar. From the aforementioned facts, using Theorem 1.3 and, if necessary,
compositions with affine isomorphisms, the following can be deduced:

 \proclaim{Corollary 1.4}  {\rm (a)} Let $\Cal L$ be a $d$-dimensional 
affine subspace of $\br^n$,   let $\mu$ be a friendly measure on $\br^d$,
and let $\vf:\br^d\to\Cal L$ be an affine isomorphism.
Then $\omega(\vf_*\mu)
= \omega(\Cal L)$.
\smallskip
 {\rm (b)} Let  $\mu$ be an absolutely decaying and Federer measure on $\br^d$, $\Cal L$ an  affine
subspace of $\br^n$,   and let
$\vf:\br^d\to\Cal L$ be a smooth map which is nondegenerate in $\Cal L$ at $\mu$-a.e.\ 
point of $\br^d$. 
Then $\omega(\vf_*\mu)
= \omega(\Cal L)$.
\endproclaim

For the special case of $\Cal L$ being extremal, part (a) was stated without proof in
\cite{KLW, \S 10.5}.


\heading{2. 
Weighted posets 
and quantitative
nondivergence}
\endheading 

In this section we  work with  mappings of weighted partially ordered
sets ({posets}) into spaces of functions on balls in a \be\ metric space.
Here is some relevant terminology.
By a {\sl weighted poset\/} we mean a partially ordered
set $\frak P$ together with a map $\eta:\frak P\to \br_+$. A linearly ordered subset of  $\frak P$
will be called a {\sl flag\/}.
We will denote by $ \ell(\frak P)$ the
{\sl length\/}
of $\frak P$, i.e.~the cardinality of a flag with   maximal cardinality. 
If $\frak F$ is a subset of $\frak P$, we let $\frak P(\frak F)$ be
the poset of elements of
$\frak P\smallsetminus \frak F$ comparable with every element of $\frak F$. Note
that one always has
$$
 \ell\big(\frak P(\frak F)\big) \le  \ell(\frak P) -  \ell(\frak F)\,.\tag 2.1
$$

We will fix a metric space $X$ and consider a weighted poset $(\frak P, \eta)$ together
with  a mapping  $\psi$ from
$\frak P$ to the space $C(B)$ of $\br$-valued
      continuous
functions on some subset $B$ of
$X$, which we will denote by
$s\mapsto \psi_s$. Given such a mapping and a positive number $\vre$, 
we will say that a
point
$z\in B$  is {\sl  $\vre$-marked 
relative to $\frak P$\/} if
there exists a flag
$\frak F_z\subset\frak P$ such that

\roster
\item"(M1)" $\vre\eta(s)\le |\psi_s(z)| \le \eta(s) \quad \forall\,s\in\frak F_z$;

\item"(M2)" $|\psi_s(z)| \ge \eta(s)\ \, \qquad\qquad\forall\,s\in \frak P(\frak F_z)$.
\endroster
We will denote the set of all such points by $\p(\vre,
\frak P)$.
When it does not cause confusion, we will omit the
reference to either $\frak P$,  $\eta$ or $\vre$, and will simply say that
$z$  is
$\vre$-marked, or marked relative to $\frak P$.

         \proclaim{Theorem 2.1}  Let 
${k},N\in
\bz_+$ and
$C,\alpha,D  > 0$.
Suppose that we are given  an 
$N$-\be\ metric space $X$,
a weighted poset $(\frak P, \eta)$, a ball $B = B(x,r)$ in
$X$,  a measure 
     $\mu$ which is $D$-Federer  on $\tilde B\df B\big(x,3^mr\big)$, and    a mapping
$\psi:\frak P\to C(\tilde B)$, $s\mapsto \psi_s$,   such that the
following holds:

\roster
\item"(A0)" $ \ell(\frak P) \le {k}$;
\item"(A1)" $\forall\,s\in \frak P\,,\quad \psi_s$ is \cag\ on $\tilde B$ with
respect to
$\mu$;
\item"(A2)" $\forall\,s\in \frak P\,,\quad\|\psi_s\|_{\mu,B} \ge\eta(s)
$;
\item"(A3)"  
$\forall\,y\in \tilde
B \,\cap\,\supp\mu,\quad\#\{s\in \frak P\bigm|
|\psi_s(y)| < \eta(s)\} < \infty$.
\endroster
Then
$
\forall\,\vre > 0$ one has
$$
\mu\big(B\smallsetminus \p(\vre,
{\frak P})\big) \le {k}C
\big(ND^2\big)^{k}
\vre^\alpha
\mu(B)\,.
$$
\endproclaim

The proof given below is a weighted modification of the argument 
from  \cite{KT, \S 5}, which, in its turn, generalizes \cite{KM, \S 4}. In a sense,
this modification allows one to use the full strength of the construction originally 
introduced by Margulis \cite{Mr1}, obtaining what may be considered as the optimal
result (see the discussion after the proof of Theorem 2.2).

\demo{Proof}  We proceed by induction on ${k}$. If
${k}=0$, the poset $\frak P$ is empty, and for any
$z\in B$  one can take
$\frak F_z = \vrn$ and check that (M1) and (M2) are satisfied for all
positive $\vre$; thus all points of $B$ are marked.  Now
take
${k}\ge 1$ and  suppose that the claim is proved for all smaller values of ${k}$.

Fix $C,\alpha,\rho,  \frak P, B = B(x,r)$ and $\psi$  as in the
formulation of
the theorem.   For any  $y\in B \,\cap\,\supp\mu$ define
$$
H(y)\df\{s\in \frak P\bigm| |\psi_s(y)| < \eta(s) \}\,;
$$
this is a finite subset of $\frak P$ in view of (A3). If
$H(y)$ is empty, $y$ is clearly $\vre$-marked  for any
positive
$\vre$: indeed, since  $|\psi_s(y)| \ge
\eta(s)
$ for all $s\in \frak P$, one can again take $\frak F_y$ to be the
empty set and
check that (M1) and (M2) are satisfied. Thus one only needs to consider
points
$y$ from the set
$$
E\df \{y\in B \,\cap\,\supp\mu \mid H(y)\ne\vrn\} = \{y\in B
\,\cap\,\supp\mu\mid
\exists\,s\in \frak P
\text{ with }|\psi_s(y)| < \eta(s) \}\,.
$$

Take $y\in E$ and $s\in H(y)$, and define
$$
r_{s,y}\df\sup\{t > 0\bigm| \|\psi_s\|_{\mu,B(y,t)} < \eta(s) \}\,.\tag 2.2
$$
         It follows from the continuity of functions $\psi_s$ that
for small enough positive $t$ one has $\|\psi_s\|_{\mu,B(y,t)} < \eta(s) $,
hence
$r_{s,y} > 0$. Denote $B(y,r_{s,y})$ by
$B_{s,y}$. From (A2) it is clear
that $B_{s,y}$ does not contain $B$; therefore one has
$r_{s,y} < 2r$.
Note also that (2.2) immediately
implies that
$$\|\psi_s\|_{\mu,B_{s,y}}\le\eta(s)\,.\tag 2.3$$


Now for any $y\in E$ choose an element $s_y$ of $ H(y)$ such that
$r_{s_y,y}
\ge   r_{s,y}$ for all $s\in H(y)$ (this can be done since $H(y)$ is finite).
For brevity let us denote $r_{s_y,y}$ by $r_y$ and $B_{s_y,y}=
\bigcup_{s\in H(y)}B_{s,y}$ by $B_y$. Also let us denote the poset
$\frak P(\{s_y\})$ by $\frak P_y$.

      \medskip
The next claim allows one to show  a
point
$z\in B_y$ to be
marked relative to $\frak P$ once it is
marked relative to $\frak P_y$. 
Namely, fix  $\vre > 0$ and $y\in E$, and take $z$ such that
$$z\in B_y\,\cap\,\supp\,\mu\,\cap \,\p(\vre,
{\frak P_y})\quad\text{and}\quad|\psi_{s_y}(z)| \ge\vre\eta(s_y)\,.\tag 2.4
$$
Then we claim that $z$ belongs to $\p(\vre, {\frak P})$; equivalently,
$$
(B_y\cap\,\supp\mu)\ssm\p(\vre,
{\frak P}) \subset
\big(B_y\ssm\p(\vre,
{\frak P_y}) \big) \cup 
\big\{ x\in B_y\bigm|
|\psi_{s_y}(x)| < \vre\eta(s_y)\}
\,.\tag 2.5
$$

Indeed, take $z$ as in (2.4); by definition of $\p(\vre,
{\frak P_y})$, there exists
a flag $\frak F_{y,z}\subset\frak P_y$ such that
$$
\vre\eta(s)\le |\psi_s(z)| \le \eta(s)\quad \forall\,s\in\frak F_{y,z}\tag 2.6
$$
and
$$
|\psi_s(z)| \ge \eta(s)\quad
\forall\,s\in \frak P_y(\frak F_{y,z})\,.\tag 2.7
$$
Put $\frak F_{z}\df\frak F_{y,z}\cup\{s_y\}$. Then
$\frak P(\frak F_{z}) = \frak P_y(\frak F_{y,z})$, so (M2)
immediately follows from
(2.7). Property  (M1) for $s  \ne s_y$ is given by (2.6), and for $s = s_y$
 by (2.4)  and
 (2.3). Thus (2.5) is proved.

      \medskip
Note that for any $y\in E$ one clearly has $r_{y} <
2r$, which in particular implies that $B_y\subset B(x,3r)$. We are  going
to  fix some
$r'_y$ strictly between
$r_y$ and
$\min(2r, 3r_{y})$, and denote $B(y,r'_{y})$ by $B'_y$. Clearly one
has
$$
\|\psi_s\|_{\mu,B'_{y}}  \ge \eta(s)
\quad\text{for any }y\in E\text{ and }s\in \frak P\,.\tag 2.8
$$
(Indeed, the definition of $r_y$ and
        (2.2) imply the above inequality for any  $s\in H(y)$,
and it obviously holds if $s\notin H(y)$.)

        Now observe that
$\frak P_y$,  $B'_y$ and $\tilde B'_y\df B(y,3^{{k}-1}r'_y)$ satisfy
properties
\roster
\item"$\bullet$" (A0) with ${k}$ replaced by ${k}-1\qquad\ \ $ --- $\qquad$in view of
(2.1);
\item"$\bullet$" (A2)
\hskip 1.9in ---
$\qquad$in view of (2.8);
\item"$\bullet$" (A1) and (A3) \hskip 1.31in ---
$\qquad$since
\endroster
$$
\tilde B'_y =  B(y,3^{{k}-1}r'_y)\subset B(x,3^{{k}-1}r'_y + r)
\subset B\big(x,(2\cdot 3^{{k}-1} + 1)r\big)\subset B(x,3^{k}r) = \tilde B\,.
$$

Therefore one has
$$
\aligned
\mu\big(B_y\ssm\p(\vre,
{\frak P_y}) \big) &\le \mu\big(B'_y\ssm\p(\vre,
{\frak P_y}) \big)\le({k}-1)C(ND^2)^{{k}-1} \vre^\alpha
\mu(B'_y)\\ &\le D(k-1)C(ND^2)^{{k}-1} \vre^\alpha
\mu(B_y)
\endaligned\tag 2.9
$$
         by the induction assumption and the Federer property of $\mu$. On the
other hand, in view of
$\psi_{s_y}$ being \cag\ on $\tilde B\supset B'_y$ with respect to $\mu$,  one can write
$$
\aligned
\mu\left(
\big\{ x\in B_y\bigm|
|\psi_{s_y}(x)| < \vre\eta(s_y)\}
\right) \le
\mu\left(
\big\{ x\in B'_y\bigm|
|\psi_{s_y}(x)| < \vre\eta(s_y)\}
\right) \\ \le
\ C\left(\frac{\vre\eta(s_y)}{\|\psi_{s_y}\|_{\mu,B'_y}}\right)^\alpha  \mu(B'_y)
\underset{(2.8)}\to\le C  \vre^\alpha
\mu(B'_y)\underset{\text{Federer}}\to \le CD \vre^\alpha
\mu(B_y)\,.\quad\endaligned\tag 2.10
$$
Now recall that we need to estimate the measure of
$E\ssm \p(\vre,
{\frak P})$. For any $y\in E$, in view of  (2.5), (2.9) and (2.10) one
has
$$
\aligned
\mu\big(B_y\ssm\p(\vre,
{\frak P})\big) &\le
C\left(({k}-1)N^{{k}-1}D^{2{k}-1} +
D \right)\vre^\alpha  \mu(B_y)%
\\ &\le
{k}CN^{{k}-1}D^{2{k}-1}\vre^\alpha  \mu(B_y)\,.
\endaligned
\tag 2.11
$$

Finally, consider the covering $\{B_y\mid y\in E\}$ of $E$, choose a
countable subset $Y$ of $E$ such that the multiplicity of the subcovering
$\{B_y\mid {y\in Y}\}$ is   at most $N$, and write
$$
\sum_{y\in Y}\mu(B_y) \le N\mu\big(\bigcup_{y\in Y} B_y\big)
\le N\mu\big(B(x,3r)\big) \le ND\mu(B)\,.\tag 2.12
$$
Therefore the measure of $E\ssm \p(\vre,
{\frak P})$ is bounded from above by
$$
\split
\sum_{y\in Y}\mu
\big(B_y\ssm\p(\vre,
{\frak P})\big)
&\underset{(2.11)}\to\le  {k}CN^{{k}-1}D^{2{k}-1}\vre^\alpha
\sum_{y\in Y}\mu(B_y)\\ &\underset{(2.12)}\to\le  {k}C\big(ND^2\big)^{{k}}\vre^\alpha
\mu(B)\,,
\endsplit
$$
and the theorem is proven. \qed\enddemo



We now apply Theorem 2.1 to the (appropriately weighted) poset 
$$
\frak P_k \df \{\text{nonzero  primitive subgroups of }\bz^k\}\,,
$$
 with the inclusion relation  (recall that  a
discrete subgroup $\Gamma\subset\bz^k$ is called {\sl primitive\/} 
if $\Gamma =
\Gamma_\br\cap \bz^k$), 
and prove a general version of Theorem  0.2.

 \proclaim{Theorem 2.2} Let $k,N\in\bn$ and $C, D, \alpha, \rho > 0$,   and
suppose we are given an  $N$-\be\ metric space $X$, a ball $B = B(x_0,r_0)\subset
X$, a measure $\mu$ which is $D$-Federer  on $\tilde B \df B\big(x_0,3^kr_0\big)$,
 and
a map
$h:\tilde B \to \GL_k(\br)$.  
 Assume that the following two conditions hold:
\roster
\item"{\rm [2.2-i]}" $
\forall\,\Gamma\in\frak P_k , \text{ the function
$x\mapsto\|h(x)\Gamma\|$ is \cag\ on $\tilde B$ w.r.t.\ $\mu$;}$
\item"{\rm [2.2-ii]}" $
\forall\,
\Gamma\in\frak P_k, \quad
\|h(\cdot)\Gamma\|_{\mu,B} \ge \rho^{\rk(\Gamma)}\,.$
\endroster
Then 
 for any  positive $ \vre \le \rho$ one has
$$
\mu\big(\{x\in B\mid 
h(x)\bz^k\notin K_\vre\}\big)
 \le 
kC \big(ND^2\big)^k \left(\frac{\vre} \rho \right)^\alpha  \mu(B)\,.
$$
\endproclaim

\demo{Proof}  We let $\frak P = \frak P_k$ and for all $ \Gamma\in\frak P$
define 
$\psi_{\sssize \Gamma}(\cdot)
\df
\|h(\cdot)\Gamma\|$  and  
$\eta(\Gamma) =  \rho^{\rk(\Gamma)}$. It is easy to verify that 
$(\frak P, \eta)$ and $\psi$ satisfy properties (A0)--(A3) of Theorem 2.1.
Indeed, (A0) holds since  any two primitive subgroups are either
incomparable or have the same rank, (A1) is given
by [2.2-i] and  (A2) by [2.2-ii]. To check  (A3) it suffices to observe  that for any $x$, since
$\bigwedge\big(h(x)\bz^k\big)$ is discrete in $\br^k$,  the number of primitive subgroups
$\Gamma$ of $\bz^k$ for which $\|h(x)\Gamma\| \le 1$ is finite. 

In view of Theorem 2.1, it remains to prove that a point $x\in B$ with 
$h(x)\bz^k\notin K_\vre$
cannot be 
$\frac\vre\rho$-marked; in other words,
$$
\p(\vre/\rho,\frak P) \subset  \big\{x\in B\bigm|
\|\vv\| \ge \vre\text{ for all }\vv\in\bz^k\nz\big\}\,.\tag 2.13
$$
Take an $\frac\vre\rho$-marked point $x\in B$, and let
$\{0\} =
\Gamma_0 \subsetneq \Gamma_1 \subsetneq\dots \subsetneq\Gamma_j = \bz^k$
be all the elements of $\frak F_x \cup \big\{\{0\},\bz^k\big\}$. Properties (M1) and (M2)
translate into:
$$
\frac\vre\rho\cdot\rho^{\rk(\Gamma_i)}\le \|h(x)\Gamma_i\| \le \rho^{\rk(\Gamma_i)} \quad
\forall\,i = 0,\dots,j-1\,,\tag 2.14
$$
and
$$
\|h(x)\Gamma\| \ge \rho^{\rk(\Gamma)}\ \, \qquad\qquad\forall\,\Gamma\in \frak P(\frak
S_x)
\,.\tag 2.15
$$
(Even though $
\Gamma_0 = \{0\}\notin\frak P$, it is clear that it also satisfies
(2.14).)

Take
any  $\vv\in\bz^k\nz$. Then there exists $i$, $1 \le i \le j$,
such that
$\vv\in \Gamma_i\ssm \Gamma_{i-1}$. Denote $(\Gamma_{i-1} +
\bz\vv)_\br\cap
\bz^k$ by $\Gamma$. Clearly it is a primitive
subgroup  of $\bz^k$ satisfying 
$\Gamma_{i-1}\subset \Gamma \subset \Gamma_i$, therefore $\Gamma \in \frak F_x \cup
\frak P(\frak F_x)$. Now one can use properties (2.14) and (2.15) to deduce that
$$
 \|h(x)\Gamma\| \ge
\min\left(\frac{\vre}
\rho \cdot\rho^{\rk(\Gamma)},\, \rho^{\rk(\Gamma)}\right)  = {\vre}
\rho^{\rk(\Gamma) - 1} = {\vre}\rho^{\rk(\Gamma_{i-1})}\,.\tag 2.16
$$
On the other hand, from the submultiplicativity of the covolume it follows that
$ \|h(x)\Gamma\|\le \|h(x)\Gamma_{i-1}\|\cdot \|h(x)\vv\|$. Thus
$$
\|h(x)\vv\| \ge \frac{\|h(x)\Gamma\|}{\|h(x)\Gamma_{i-1}\|} \underset{\text{by (2.14) and
(2.16)}}\to\ge 
\frac{{\vre}
\rho^{\rk(\Gamma_{i-1})}} { \rho^{\rk(\Gamma_{i-1})} }  =  \vre\,.
$$
This shows (2.13) and completes the proof of the theorem. \qed\enddemo

In order to better understand the difference between this theorem and
its predecessors (proved in \cite{KM, KLW, KT}), let us draw a corollary from
it. Namely, suppose that $X$, $\mu$ and $B$ are as in the above theorem, that $h$ satisfies
[2.2-i] with some $C,\alpha$ (for example, $X = \br^d$, $\mu = \lambda$ and $h$ is a
polynomial map), and that for some small positive $\vre$, the relative measure
of $x\in B$ for which $h(x)\bz^k\notin K_\vre$ is  at least $1/2$. Then Theorem 2.2
asserts that there exists a subgroup $
\Gamma$ of $\bz^k$  such that 
$\|h(\cdot)\Gamma\|_{\mu,B}< \rho^{\rk(\Gamma)}$, where $\rho = \left( 2kC
\big(ND^2\big)^k
 \right)^{1/\alpha}\vre$. Consequently, in view of Minkowski's Lemma, the whole
`trajectory' $h(B\cap\supp\mu)\bz^k$ must be contained in the complement to
$K_{\const\cdot\vre}$, with the constant depending only on $D, N, C, \alpha$ and $k$.
In other words, it must stay at a (uniformly) bounded distance from the complement to
$K_\vre$. Note that using previously known results it was only possible to conclude that
$h(B\cap\supp\mu)\bz^k$ must be outside of
$K_{\const\cdot\vre^{1/(k-1)}}$, a compact set of diameter approximately
$(k-1)$ times smaller than that of $K_\vre$.

When it comes to number-theoretic applications, the crucial advantage is that
replacing $\rho$ by $\rho^{\rk(\Gamma)}$ 
makes [2.2-ii] easier to check. This will be demonstrated in the next section,
where Theorem 2.2 will be applied to $h$ as in (0.4). 


\heading{3. An application to metric \da}
\endheading

We recall some notation and terminology introduced in the beginning of this paper,
as well as in the paper \cite{K2}.   
For $m,n\in\bn$ and $v > 0$, we denote by $\Cal W_v$ the set of \amr\  for which  there are
infinitely many
$\vq\in
\bz^n$ such that (0.5) holds for some $\vp\in\bz^m$.
The dimensionality of the matrices will be clear from the context. 
The {\sl
\de\/} $\omega(A)$ of $A$ 
 defined in the Introduction is equal to
$$
\omega(A) = \sup\{v\mid A\in \Cal W_v\}\,.
$$
Clearly $\Cal
W_{u} \subset \Cal
W_{v}$ if $u \ge v$.
We will also use the notation
$$
\Cal W_v^{\sssize +}\df \bigcup_{u > v} \Cal W_{u} = \{A\mid
\omega(A) > v\}\,.
$$
Note that the definition of $\omega(A)$, unlike that of the sets  $\Cal W_v$, does not 
depend on the choice of norms on
$\br^m$ and
$\br^n$.

Although there are many interesting and unsolved \di\ problems related to 
the space of $m\times n$ matrices, we specialize to the case $m = 1$, 
that is, consider \di\ properties
of vectors ($=$ row matrices) $\vy\in\br^n$. With some abuse of notation, 
we will view integer vectors 
$\vq\in\bz^n$ as column vectors, so that 
$\vy\vq$  stands for 
$y_1q_1 + \dots + y_nq_n$.

Now let us describe a  correspondence, dating back to \cite{Sc1} and \cite{Da3}, 
between approximation properties of vectors  $\vy\in\br^n$ and dynamics of
certain trajectories  in $\Omega_{n+1}$. Given a
row vector  $\vy\in\br^n$ one  defines 
$$
u_\vy \df \left(\matrix
1 & \vy  \\
0 & I_n
\endmatrix \right)\,,\tag 3.1
$$
and
considers
 the lattice  $u_\vy\bz^{n+1}$ in $\br^{n+1}$, 
that is, the collection of vectors  of the form 
$\left(\matrix \vy\vq + p \\  \vq \endmatrix
\right)
$,
where $p\in\bz$ and 
$\vq \in \bz^n$. Then one can read \di\ properties of
$\vy$  from the behavior of the trajectory 
$Fu_\vy\bz^{n+1}$ in the space of lattices,
where
$$
F = \{g_t\mid t\ge 0\}\,,\quad\text{with}\quad g_t = \text{\rm diag}(e^{t},
e^{-t/n},\dots,e^{-t/n})\,,\tag 3.2
$$
is  a one-parameter subsemigroup of
$\SL_{n+1}(\br)$  which  expands the first  coordinate and uniformly contracts the last $n$
coordinates of vectors in $\br^{n+1}$.  

\medskip
The  following 
elementary lemma was proved in \cite{K2}.

\proclaim{Lemma  3.1} Suppose we are given a set $E\subset \br^2$ which is
discrete and homogeneous with respect to positive integers, that is,
 $kE\subset E$ for any $k\in\bn$. Also take $a,b > 0$, $v > a/b$, and define $c$ by 
$$
c = \frac
{bv - a}{v + 1}\quad \Leftrightarrow 
\quad v = \frac{a+c}{b-c}\,.
$$ 
Then the following are equivalent:

\roster
\item"{\rm [3.1-i]}" 
there exist $(x,z)\in E$ with
arbitrarily large
$|z|$ such that 
$|x| \le |z|^{-v}$;
\item"{\rm [3.1-ii]}" there exist arbitrarily large $t > 0$ such that for some
$(x,z)\in E\nz$ one has
$$
\max\left(e^{a t} |x|,  e^{-bt} |z|\right) \le e^{-c t}\,.
$$ 
\endroster
\endproclaim

Taking $v > n$, $\vy\in\br^n$ and
$$
E = \big\{\big(\vy\vq + p,\|\vq\|\big)\bigm| (p,\vq)\in\bz^{n+1}\big\}\,,
$$
one  notices that [3.1-i] is equivalent to $\vy\in\Cal W_v$. On the other hand, 
choosing $a = 1$ and  $b = 1/n$ one sees that [3.1-ii] amounts to 
$$
 g_{t} u_\vy\bz^{n+1}\notin K_{ e^{-c t}}
\text{\ \ for an unbounded set of }t \in \br_{+}\,,\tag 3.3 
$$
  where
$$
c = \frac
{v-n}{n(v+1)}\quad \Leftrightarrow 
\quad v = \frac{n(1+c)}{1- n c }\,.\tag 3.4
$$ 
Thus, if we define $\gamma(\vy)$ 
to be the supremum of all $c$ for which (3.3) holds, or, equivalently,
$$
\gamma(\vy) = \sup\left\{c\, \left| \ 
 g_{t} u_\vy\bz^{n+1}\notin K_{ e^{-c t}}
\text{\ for infinitely many }t \in \bn\right.\right\} 
$$ (in \cite{K2} this quantity was
called  the {\sl growth exponent\/}
 of $u_\vy\bz^{n+1}$ with respect to $F$), then we have the equality 
$$
\omega(\vy) = n\frac{1+\gamma(\vy)}{1- n \gamma(\vy) }\tag 3.5
$$ 
(cf.~\cite{K2, Corollary 2.3}).

\medskip

Now let us turn to computing \de s of measures. If $\nu$ is a measure on $\br^n$ 
and $v \ge n$, one has $\omega(\nu)\le v$ iff $\nu(\Cal W_u) = 0$ for any $u > v$.
In view of the above discussion, this amounts to saying that for any $d
> c$ where $c$ is given by (3.4), 
$$
\nu\left(
\left\{\vy \left| 
 g_{t} u_\vy\bz^{n+1}\notin K_{ e^{-d t}}
\text{\ for an unbounded set of }t \in \br_{+}\right.\right\}\right) = 0\,.
$$ 
This is easily seen to be equivalent to
$$
\forall\, d > c\,,\quad\nu\left(
\left\{\vy \left| 
 g_{t} u_\vy\bz^{n+1}\notin K_{ e^{-d t}}
\text{\ for infinitely many }t \in \bn\right.\right\}\right) = 0\,,
$$ 
and, in view of the Borel-Cantelli Lemma, a condition sufficient for the latter is
$$
\sum_{t = 1}^\infty\nu\left(
\left\{\vy \left| 
 g_{t} u_\vy\bz^{n+1}\notin K_{ e^{-d t}}
\right.\right\}\right) < \infty\quad\forall\, d
> c\,.\tag 3.6 
$$ 
This is precisely where the measure estimates discussed in the previous section come in.

 \proclaim{Proposition 3.2} Let $X$ be a \be\ metric space,  $B = B(x,r)\subset
X$  a ball,  $\mu$ a measure which is $D$-Federer  on $\tilde B \df
B(x,3^{n+1}r)$ for some $D > 0$, and
$\vf$ a continuous  map from $\tilde B$ to $\br^n$. Also take $c\ge 0$ and  assume that:

\roster
\item"{\rm [3.2-i]}" $\exists\,C,\alpha > 0$ such that 
all the functions
$x\mapsto\|g_t u_{\vf(x)}\Gamma\|$,  $\Gamma\in\frak P_{n+1}$, are \cag\ on $\tilde B$
w.r.t.\ $\mu$;
\item"{\rm [3.2-ii]}" for any $\,d
> c$  there exists $T = T(d) > 0$ such that for any $t\ge T$ and
any $\Gamma\in\frak P_{n+1}$ one has
$$
\|g_t u_{\vf(\cdot)}\Gamma\|_{\mu,B}\ge e^{-{\rk(\Gamma)}d t}\,.
$$
\endroster 
Then $\omega\big(\vf_*(\mu|_B)\big) \le v$, 
where $v$ is
related to
$c$ via {\rm (3.4)}.
\endproclaim

\demo{Proof}  As was
observed  in the course of the preceding discussion, it suffices
 to
show (3.6) for $\nu = \vf_*(\mu|_B)$. We now proceed to verify that the map
$h = g_t u_{\vf}$  satisfies the assumptions of  Theorem 2.2, with $k = n+1$. 
Condition  [2.2-i] clearly coincides with [3.2-i].  For the other
condition, we take $d>c$ and $\rho = e^{-\frac{c+d}2 t}$, so that [3.2-ii] 
implies [2.2-ii] for any  $t > T(\frac{c+d}2)$. Therefore, by Theorem 2.2, 
$$
\split
\nu\left(
\left\{\vy \left| 
 g_{t} u_\vy\bz^{n+1}\notin K_{ e^{-d t}}
\right.\right\}\right)  &= \mu\big(\{x\in B\mid
h(x)\bz^{n+1}\notin K_{e^{-d t}}\}\big) \\ &\le  (n+1)C \big(ND^2\big)^{n+1} \left(e^{-d
t}e^{\frac{c+d}2 t}
\right)^\alpha  \mu(B) \\ & = \const\cdot e^{-\alpha\frac{d-c}2 t}
\endsplit
$$
for all but finitely many $t\in\bn$. This readily implies (3.6).
\qed\enddemo

\example{Remark 3.3} Note that both [3.2-i] and [3.2-ii] trivially hold for
$\Gamma$ of rank $0$ and $n+1$, and also that the validity of those conditions
for all primitive $\Gamma$ is equivalent to their validity for all subgroups. 
It will be convenient to denote by  $\Cal S_{k,{j}}$ the
set of all subgroups of $\bz^k$ of rank ${j}$. Thus, to apply 
Proposition 3.2 it will be enough to check [3.2-i,ii] for all $\Gamma\in \Cal S_{n+1,{j}}$,
${j} = 1,\dots,n$. 
\endexample

\heading{4. Computing \de s of measures}
\endheading

\comment
From this point on it will be convenient to use the following function on the space   of 
lattices in $\br^k$:
$$
\delta(\Lambda) \df \inf_{\vv\in\Lambda\nz}\|\vv\|\,,
$$
so that
$K_\vre = \big\{\Lambda  \bigm| \delta(\Lambda) \ge \vre\big\}\,.$
With this notation,
\endcomment

In this section we  use Proposition 3.2 to prove Theorem 1.3, that is,
for a given $v\ge n$, write necessary and sufficient conditions for \de s of certain
measures to be not greater than $v$.
For this we need to understand to what extent the two conditions in the above
proposition are necessary. While not much can be said  about 
the first one, it turns out that 
 assumption [3.2-ii] is in fact necessary
for the conclusion of Proposition 3.2. 
Furthermore, the consequences of  [3.2-ii] not being true
are surprisingly strong.

 \proclaim{Lemma 4.1}  Let $\mu$ be a measure on a set $B$, take
$c,v > 0$  related via  {\rm (3.4)}, and let
$\vf$ be a 
 map from $
B$  to $\br^n$ such that {\rm [3.2-ii]} does not hold. Then
$\vf(B\cap\,\supp\,\mu)\subset \Cal W_u$  for some $u > v$. 
\endproclaim

We remark that the negation of the conclusion of Proposition 3.2 is much weaker:
it simply
amounts to saying that for some $u > v$ the set $\{x\in
B\mid\vf(x)\in\Cal W_u\}$ has positive measure.

\demo{Proof} The assumption of the lemma says that 
there exists   $d > c$ such that 
one has
$$\sup_{x \in B\cap\,\supp\,\mu}\|g_t u_{\vf(x)}\Gamma\| < e^{-{\rk(\Gamma)}d t}$$
for arbitrarily large $t$ (and $\Gamma\subset \bz^{n+1}$ dependent on $t$). 
In other words, there exists $1\le {j} \le n$, a sequence ${t_i}\to\infty$ and a sequence
of subgroups
$\Gamma_i \in \Cal S_{n+1,{j}}$ such that for any $x \in B\cap\,\supp\,\mu$ one has
$\|g_{t_i} u_{\vf(x)}\Gamma_i\| < e^{-{{j}}d {t_i}}$. But
in
view of Minkowski's Lemma, this means that for any $i$ and any $x \in B\cap\,\supp\,\mu$
there is a nonzero vector
$\vv\in g_{t_i} u_{\vf(x)}\Gamma_i$ with $\|\vv\|  < 2^{j} e^{-d {t_i}}$.
Hence 
$$
 g_{t_i} u_{\vf(x)}\bz^{n+1}\notin K_{2^{j} e^{-d t_i}}
\,,$$ 
which implies that 
$\gamma\big(\vf(x)\big) \ge d$, finishing the proof in view
of 
(3.5). 
\qed\enddemo

It is worthwhile to point out that it is precisely the above argument that requires
a strengthening of the quantitative nondivergence obtained in \S 2. Previously available
techniques could only yield $\gamma\big(\vf(x)\big) \ge d/n$, which was enough in the case
$v = n$ $\Leftrightarrow$ $c = 0$, but not in the general case.

\medskip

To write [3.2-ii] in a more convenient form, instead  of 
discrete subgroups $\Gamma$ of $\br^k$ we will work with elements $\vw\in \bigwedge(\br^k)$ representing them, saying that $\vw$ {\sl represents\/} $\Gamma\ne\{0\}$ if 
$
\vw = \vv_{1}\wedge\dots\wedge \vv_{j}$ where $\vv_{1},\dots, \vv_{j}$ is a basis of $\Gamma$
as a $\bz$-module.
Clearly $\vw$ representing $\Gamma$ is defined up to a sign, hence $\|\Gamma\| = \|\vw\|$. 
Now let us reproduce a computation (first done in
\cite{KM}) of  coordinates of $g_t u_{\vf(x)}\Gamma$ with respect to the standard
basis of $\bigwedge(\br^{n+1})$. For the rest of this section let us denote $\br^{n+1}$
by $V$, its standard 
 basis  by 
$\ve_0,\ve_1,\dots,\ve_n$,
and
the space spanned by $\ve_1,\dots,\ve_n$ by $V_0$. For
$I =
\{i_1,\dots,i_{j}\}\subset \{0,\dots,n\}$, $i_1 < i_2 < \dots < i_{j}$,  let 
$$\ve_{\sssize I} \df \tsize
\ve_{i_1}\wedge\dots\wedge \ve_{i_{j}}\in \bigwedge^{j}(V)\,,$$
 with the convention
$\ve_\vrn = 1$; then $\{\ve_{\sssize I}\mid \# I = {j}\}$ is a basis of
$\bigwedge^{j}(V)$, and we extend the Euclidean structure $\langle\cdot,\cdot\rangle$ from $V$
to its exterior powers so that this basis becomes orthonormal.

Since the action of
$u_{\vy}$ leaves $\ve_0$ invariant and sends $\ve_{i}$, $i > 0$,  to $\ve_{i} + y_i
\ve_0$, one can write\footnote{The quantity  $\langle\ve_i
\wedge\ve_{\sssize I\ssm \{i\}},\ve_{\sssize I}\rangle$ in (4.1), which
is equal to $1$ (resp., $-1$) if the number of elements of $I$ strictly between $0$ and $i$ 
is even (resp., odd), was denoted by $(-1)^{\ell(I,i)}$ in \cite{K2}.}
$$
u_{\vy}\ve_{\sssize I} =  \cases &\ve_{\sssize I} \hskip 2.05in \text{ if }0\in I\,,\\
&\dsize\ve_{\sssize I} + \sum_{i\in I} \langle\ve_i\wedge\ve_{\sssize I\ssm \{i\}},\ve_{\sssize I}\rangle
 \, y_i\,
\ve_{\sssize I \cup \{0\}\ssm\{i\}}\ \text{ otherwise}\,.\endcases\tag 4.1
$$
Therefore $\vw\in \bigwedge^{j}(V)$ is sent to
$$
u_{\vy}\vw  = \sum \Sb I\subset\{1,\dots,n\}
\endSb
\langle\ve_{\sssize I},\vw\rangle \,\ve_{\sssize I} + 
\sum \Sb J\subset\{1,\dots,n\}
\endSb
\left( \sum_{i= 0}^n \langle\ve_i\wedge\ve_{\sssize J},\vw\rangle\, y_i\right)
\ve_{\sssize \{0\}\cup J} \tag 4.2
$$
(where it is understood that $\# I = j$, $\# J = j-1$, and $y_0 = 1$). Note that the first sum in (4.2) is the image of $\vw$ under the orthogonal
projection  from $\bigwedge^{j}(V)$ onto   $\bigwedge^{j}(V_0)$, which (the projection)
we will denote by $\pi$. On the other hand, each term in the second sum is orthogonal to $\bigwedge(V_0)$. 
To simplify (4.2), let us define a linear map $\vc$ from $\bigwedge^j(V)$
to $\big(\bigwedge^{j-1}(V_0)\big)^{n+1}$  by setting the $i$th component
of $\vc(\vw)$, $i = 0,1,\dots,n$, equal to
$$
\vc(\vw)_i \df \,
\sum \Sb J\subset\{1,\dots,n\}\\ \# J = j-1\endSb
\langle\ve_i\wedge\ve_{\sssize J},\vw\rangle\,\ve_{\sssize J}\,.\tag 4.3
$$
For example, the choice $i = 0$ gives $
\vc(\vw)_0 = \sum \Sb I\ni 0
\endSb
\langle\ve_{\sssize I},\vw\rangle \,\ve_{\sssize I\ssm \{0\}}$, so that
$$
\ve_0\wedge\vc(\vw)_0  = \vw - \pi(\vw)\,.\tag 4.4$$ 
In particular, the kernel of $\vw\mapsto \vc(\vw)_0$ is $\bigwedge^j(V_0)$, that is,
the orthogonal complement to $\ve_0\wedge\bigwedge^{j-1}(V_0)$. Likewise, for every $i$ 
the kernel of $\vw\mapsto \vc(\vw)_i$ is 
the orthogonal complement to $\,\ve_i\wedge\bigwedge^{j-1}(V_0)$, which implies that the kernel of 
$\vc$ is trivial.

With this notation,
(4.2) can be rewritten as $$
u_{\vy}\vw =  \pi(\vw) +  \ve_0\wedge\tilde \vy\vc(\vw) \,,$$ 
where $\tilde \vy$ stands for the row vector
$(1,y_1,\dots,y_n)$, and the product in $\tilde \vy\vc(\vw)$ is the formal matrix multiplication
(of the row vector $\tilde \vy\in\br^{n+1}$ and the column vector $\vc(\vw)\in\big(\bigwedge(V)
\big)^{n+1}$).  In other words, 
$$u_{\vy}\vw =  \pi(\vw) +  \ve_0\wedge\sum_{i = 0}^{n}y_i\vc(\vw)_i
 \un{(4.4)}= \vw + \ve_0\wedge\sum_{i = 1}^{n}y_i\vc(\vw)_i\,.$$ 
Since  $\bigwedge^j(V_0)$ and   its orthogonal complement are eigenspaces of $g_t$  
with eigenvalues  $e^{-\frac {j}n t}$ and $e^{\frac {n+1-{j}}n t}$, respectively, 
one can write
$$
g_tu_{\vy}\vw =  e^{-\frac {j}n t}
\pi(\vw) + e^{\frac {n+1-{j}}n t}  \ \ve_0\wedge\tilde \vy\vc(\vw) \,;\tag 4.5
$$
%
thus, up to a uniform
constant,
$$
\|g_t u_{\vf(\cdot)}\Gamma\|_{\mu,B} = \max\left(e^{\frac {n+1-{j}}n t} 
\|\tilde \vf(\cdot)
\vc(\vw)\|_{\mu,B}\ , e^{-\frac {j}n t} \|\pi(\vw)\| \right) \,,\tag 4.6
$$
where $\Gamma\in\Cal S_{n+1,{j}}$ is represented by $\vw$, and 
we use the notation $\tilde
\vf
\df (1,f_1,\dots,f_n)$. 

\medskip

Observe that from the above one can already extract a nice lower bound for (4.6)
whenever the restrictions of $1,f_1,\dots,f_n$ to $B\cap \supp\,\mu$ are linearly
independent. Indeed, then the correspondence $\vv\mapsto\|\tilde \vf(\cdot)
\vv\|_{\mu,B}$ yields a norm on $\big(\bigwedge(V)\big)^{n+1}$, which is obviously equivalent 
to the standard (Euclidean) norm. Thus, up to a  (uniform in $\vw$) constant, the
expression (4.6) is bounded from below by $
\|\vc(\vw)\|$. Since  $\vc$ has trivial kernel and maps $\bigwedge(V_\bz)$
into $\big(\bigwedge(V_\bz)\big)^{n+1}$,  one has $\|\vc(\vw)\|\ge 1$  
for any nonzero $\vw\in
\bigwedge(V_\bz)$.
This  argument appears in  \cite{KM} and, in a dual form, in \cite{KLW}.

In general the  desired lower bound is affected by the linear dependence relations
between the components of 
$\tilde
\vf$. Namely, denote by $\Cal F_{\mu,B}$  the 
$\br$-linear span of the restrictions of $1,f_1,\dots,f_n$ to $B\cap\,\supp\,\mu$, 
denote its  dimension by $s + 1$, and choose  functions
$g_1,\dots,g_{s}: B\cap \supp\,\mu\to \br$ such that
$1,g_1,\dots,g_{s}$ form a basis of $\Cal F_{\mu,B}$.     This choice  defines 
a matrix $$R = (r_{i,j})\Sb i = 0,\dots, s\\j = 0,\dots, n\endSb\in M_{{s}+1,n+1}$$ formed by
coefficients in the expansion of
$1,f_{1},\dots,f_n$ as linear combinations of $1,g_1,\dots,g_{s}$. 
In other words, with the notation 
$\tilde \vg \df (1,g_1,\dots,g_{s})$,
one has 
$$ 
\tilde \vf(x) = \tilde
\vg(x)  R\quad\forall\,x\in B\cap \supp\,\mu\,.\tag 4.7
$$
Moreover, since the first components of $\tilde \vf$ and $\tilde \vg$ are equal to $1$,
the elements in the first column of $R$ are
$$
r_{i,0} = \cases 1 \quad\text{ if }i = 0\,,\\ 0 \quad\text{ otherwise.}\endcases\tag 4.8
$$

In view of (4.7), 
$\|\tilde \vf(\cdot)
\vc(\vw)\|_{\mu,B}$ can be replaced by $\|\tilde \vg(\cdot)
R\vc(\vw)\|_{\mu,B}$, and the latter, in view of linear independence
of the components of 
$\tilde
\vg$, simply by the norms of vectors $R\vc(\vw)$.
 Summarizing the discussion, we see that {\rm [3.2-ii]} is  equivalent to
$$
\aligned
\forall\,d
> c\ \ \exists\, T  > 0\ \text{  such that   }\ \forall \,t\ge T\,, \ \forall\,{j} =
1,\dots,n
 \text{ and
} \forall\,\vw \in \Cal S_{n+1,{j}}\\ \text{one has}\quad
\max\left(e^{\frac {n+1-{j}}n t}  \|R
\vc(\vw)\|, e^{-\frac {j}n t}  \|\pi(\vw)\| \right)\ge e^{-{j}d t}\,,\qquad
\endaligned\tag 4.9
$$
where we have identified  $\Cal S_{n+1,{j}}$ with the set of elements of
$\bigwedge^{j}(\br^{n+1})$ representing $\Gamma\in\Cal S_{n+1,{j}}$.

Here is another way to understand the above condition. Let $\vr_i = (r_{i,0},\dots,r_{i,n})$ stand for 
the $i$th row of $R$, $i = 0,\dots s$. Then, using (4.3), one can write the $i$th component 
of $R
\vc(\vw)$
in the form
$$
\aligned
[R\vc(\vw)]_i  
&
= \sum_{k=0}^n r_{i,k}
\sum \Sb J
\endSb
\langle\ve_k\wedge\ve_{\sssize J},\vw\rangle\,\ve_{\sssize J} 
\\&
= \sum \Sb J
\endSb
\big\langle \sum_{k=0}^nr_{i,k}\ve_k\wedge\ve_{\sssize J},\vw\big\rangle\,\ve_{\sssize J}
 = \sum \Sb J
\endSb
\langle \vr_i\wedge\ve_{\sssize J},\vw\rangle\,\ve_{\sssize J}\,,
\endaligned\tag 4.10
$$
and therefore, up to a uniform constant induced by replacing the Euclidean norm with the sup norm,
$$
\|R
\vc(\vw)\| = \max_{i = 0,\dots,s} \max\Sb J
\subset\{1,\dots,n\}\\ \# J = j-1
\endSb|\langle \vr_i\wedge\ve_{\sssize J},\vw\rangle|\,.\tag 4.11
$$

At this point it becomes useful to recall Lemma 3.1. Namely,  for each ${j} =
1,\dots,n$  consider
$$E = \big\{\big( \|R
\vc(\vw)\|,\|\pi(\vw)\|\big) \bigm|\vw \in \Cal S_{n+1,{j}}\big\}\,.$$
It is clearly homogeneous with respect to positive integers, and 
the fact that it is discrete is easy from (4.11) and (4.8).
Take $a = \frac {n+1-{j}}n $ and $b = \frac {j}n$, and recall that 
in the beginning we fixed $v\ge n$ and chose $c =\frac
{v-n}{n(v+1)}$. Thus 
 (4.9) amounts to saying that for any $c >  c_0 \df {j}\frac
{v-n}{n(v+1)}$ condition [3.1-ii] does not hold.
By Lemma 3.1, this is equivalent to saying that [3.1-i] does not hold with 
$v$ replaced by any real number greater than
$$
\frac{a + c_0}{b - c_0} = \frac{\frac {n+1-{j}}n + {j}\frac
{v-n}{n(v+1)}}{\frac {j}n - {j}\frac
{v-n}{n(v+1)}} = \frac {v+1-{j}}{j}\,.
$$
Therefore (4.9) becomes equivalent to 
$$
\aligned\forall\,{j} =
1,\dots,n,\ \forall\,u > \tfrac {v+1-{j}}{j} \text{ and }\forall\, \vw \in \Cal S_{n+1,{j}}\\
\text{   with large enough }
 \|\pi(\vw)\|\,, \ 
\text{  one has }\ \|R\vc(\vw)\| > \|\pi(&\vw)\|^{-u}\,.
\endaligned\tag 4.12
$$

As a result, we managed to get rid of an auxiliary variable $t$ in (4.9) and
found a way to directly involve $v$, rather than  relate it to $c$ via  (3.4).

\medskip

Note that the only way the ball $B$, the measure $\mu$  and the map $\vf$ enter the above
conditions is via the matrix $R$, which depends on both $B\cap 
\supp\mu$ and  $\vf$  and is not
uniquely  determined -- but another choice of $R$ would clearly yield a condition
equivalent to (4.9)$\Leftrightarrow$(4.12).
Here is another useful way to describe  $R$. Let $$\Cal L = \langle\vf(B\cap 
\supp\mu)\rangle_a\,,\tag 4.13$$ 
put $s = \dim(\Cal L)$, and suppose 
$$
\vh:\br^{s}\to\Cal L\text{ is an affine isomorphism, and}\quad \tilde
\vh(\vx) =  \tilde \vx
R\,,\ \vx\in\br^{s}\,,\tag 4.14
$$
where as usual we have $\tilde
\vh \df (1,h_1,\dots,h_n)$ and $\tilde
\vx \df (1,x_1,\dots,x_s)$. Then it is clear that $R$ and
$\vg\df\vh^{-1}\circ
\vf$ satisfy (4.7), and that  $1,g_1,\dots,g_s$ generate $\Cal F_{\mu,B }$ and
are linearly independent over $\br$. This way, condition 
(4.12)$\Leftrightarrow$(4.9)$\Leftrightarrow$[3.2-ii] becomes a property of 
the `enveloping subspace' 
$\langle\vf(B\cap 
\supp\mu)\rangle_a$; in particular,   $R$ can be  chosen  uniformly for all  measures
$\mu$, balls $B$ intersecting
$\supp\,\mu$ and maps $\vf$ as long as (4.13) holds.

\comment
Furthermore, in this case 
 as above can easily
exhibited if one knows that the image of 
$\vf$ lies   in some affine subspace
$\Cal L$ of
$\br^n$ in which $\vf$ is $\mu$-nonplanar. Indeed, the latter means that for any ball
$B$ intersecting
$\supp\,\mu$, the image of $B\cap\supp\,\mu$ is not contained in any proper subspace
of $\Cal L$.

Specializing to the class of good maps, we can now formulate the  desired
`if-and-only-if' statement.
\endcomment

We are now ready for the main result of the section.

 \proclaim{Theorem 4.2}  Let  $\mu$ be a Federer measure on a \be\ metric space $X$,
$\Cal L$ an affine subspace of
$\br^n$, and let $\vf:X\to\Cal L$ be a continuous map such that $(\vf,\mu)$ is good
and
nonplanar in $\Cal L$. 
Then the following are equivalent for  $v \ge n$:

\roster

\item"{\rm [4.2-i]}"  
$\{x\in
\supp\,\mu\mid\vf(x)\notin\Cal W_u\}$ 
is nonempty for any $u > v$;
\item"{\rm [4.2-ii]}" $\omega(\vf_*\mu) \le v$ (in other words, each of the 
sets in {\rm [4.2-i]} has full measure);
\item"{\rm [4.2-iii]}" 
{\rm (4.12)} holds for some ($\Leftrightarrow$ for
any)
$R$ satisfying {\rm (4.14)}.
\endroster 
\endproclaim

This was proved in \cite{K2} for $v = n$.

\demo{Proof} Obviously [4.2-ii]$\Rightarrow$[4.2-i].
Assuming [4.2-iii] and using the facts that $\mu$ is Federer and $(\vf,\mu)$ is good, 
one can conclude that $\mu$-a.e.~$x\in X$ 
has a neighborhood $V$  such that $\mu$ is
\cag\ and $D$-Federer on $V$ for some 
$C,D,\alpha > 0$.
Choose  a ball $B = B(x,r)$ of positive measure such that the
dilated ball
$\tilde B = B(x,3^{n+1}r)$ is contained in $V$, and note that (4.13) holds in view of (1.3). 
We have seen in (4.5) that for any
$\vw$, each of the coordinates of
$
g_tu_{\vf}\vw$ is expressed as a linear combination of functions
$1,f_1,\dots,f_n$. Therefore, in view of Lemma 1.1,  [3.2-i] will hold  (perhaps with a
different constant $C$), and  [3.2-ii] was postulated in the disguise of (4.12), as shown
by the discussion preceding the statement of the theorem. Thus Proposition 3.2
applies, and  [4.2-ii] follows.

On the other hand, saying that [4.2-iii] does not hold, in view
of the same discussion, amounts to assuming that no ball $B$ intersecting
$\supp\,\mu$ satisfies [3.2-ii]. 
From Lemma 4.1 it then follows that $\vf(B\cap\,\supp\,\mu)\subset \Cal W_u$  for some
$u > v$, 
contradicting [4.2-i]. 
 \qed\enddemo

An especially remarkable feature of the above theorem is that whenever conditions 
[4.2-i,ii] hold for {\it some\/} $X$,  $\mu$ and $\vf$ satisfying the assumptions of
the theorem, they hold for {\it any\/}  $X$,  $\mu$, $\vf$ satisfying those
assumptions. Indeed, condition [4.2-iii] equivalent to them has no reference to $\vf$, 
$\mu$ or $X$, only to $\Cal L$. In particular, one can make the most natural 
choice of  $X$,  $\mu$, $\vf$; that is, put $X = \br^s$, $\mu = \lambda$ and $\vf = \vh$
as in (4.14), thus establishing Theorems 1.3 and 0.3, and furthermore producing
a formula for the \de\ of $\Cal L$. Namely, one has

 \proclaim{Corollary 4.3}  Let $\Cal L$ be an ${s}$-dimensional affine subspace 
 of $\br^n$. Then $$
\omega(\Cal L) = \max\big(\,n,\ \sup\{v\mid \text{\rm (4.12) holds for $R$ as in
(4.14)}\}\big)\,.
$$ 
\endproclaim

This will be made more explicit  in the next section.

\heading{5. Higher \de s and Theorem 0.4}
\endheading 

In this section we start by fixing a  parametrization (0.6) for an  
${s}$-dimensional affine subspace  $\Cal L$  of $\br^n$.
This amounts to taking $R$ of the form $$R = R_{\sssize A} \df\pmatrix I_{s+1} &
A\endpmatrix$$ where $A\in M_{s+1,n-s}$. In order to restate condition (4.12) in terms of
$A$, let us denote by $\va_i = (a_{i,s+1},\dots,a_{i,n})$ the $i$th row of $A$, $i = 0,\dots,s$,
and identify it with an element of $V$ by putting $\va_i = \sum_{k = s+1}^n a_{i,k}\ve_k$.
Then, in view of  (4.11), we have $$
\|R_{\sssize A}
\vc(\vw)\| = \max_{i = 0,\dots,s} \max\Sb J
\subset\{1,\dots,n\}\\ \# J = j-1
\endSb\big|\big\langle (\ve_i + \va_i)\wedge\ve_{\sssize J},\vw\big\rangle\big|\,.\tag 5.1
$$
Corollary 4.3  asserts that for $v\ge n$, the \de\ of $\Cal L$ is not greater than $v$ 
if and only if for all
${j} =
1,\dots,n$, $u > \tfrac {v+1-{j}}{j}$ and $\vw \in \Cal S_{n+1,{j}}$
  with large enough $\|\pi(\vw)\|$, the expression in (5.1) 
exceeds $\|\pi(\vw)\|^{-u}$.

This condition however can be considerably simplified. 
 Namely, put 
$$ V_{\bullet} \df\br\ve_{s+1}\oplus\cdots\oplus\br\ve_{n}\,,$$ and denote by 
$\pi_{\bullet}$ 
the projection of $\bigwedge(V)$ to $\bigwedge(V_{\bullet})$. We will also be using 
the following notation: $x \ll y$ will stand for $x < Cy$, where $C$ depends only 
on the matrix $A$ and  not   on $\vw$.

 \proclaim{Lemma 5.1} Suppose that $\|R_{\sssize A}
\vc(\vw)\|$ is less than $1$ for some $\vw\in \Cal S_{n+1,{j}}$.
Then  $\|\vw\|\ll1 + \|\pi_{\bullet}(\vw)\|$.  
\endproclaim

\demo{Proof} Let us take $I\subset \{0,\dots,n\}$ of size $j$ and  prove that  the absolute value of 
$\langle \ve_{\sssize I},\vw\rangle$ is bounded from above by a uniform constant times $\pi_{\bullet}(\vw)$.
Denote by $k$ the smallest element
of $I$. The claim  is trivial if $k >s$. Otherwise,  using    (5.1),
 one can write
$$|\langle \ve_{\sssize I},\vw\rangle| < 1 + \big|\langle  \va_k
\wedge\ve_{\sssize I\ssm\{k\}},\vw\rangle\big|$$
and observe that, since $\va_i\subset V_{\bullet}$ for each $i$, the 
 right hand side is not greater than $1 + \max_i\|\va_i\|\cdot\max_{J\subset\{k+1,\dots,n\}}
|\langle \ve_{\sssize J},\vw\rangle|$. The same argument can be applied
to each of the components $\langle \ve_{\sssize J},\vw\rangle$, and after 
no more than $s$ additional steps the process will terminate.   \qed\enddemo

In particular, the lemma forces  $\|R_{\sssize A}
\vc(\vw)\|$ to be not less than $1$ for any  
 $\vw\in \Cal S_{n+1,{j}}$ with $j > n-s$   and large enough $\pi(\vw)$ (this was checked in \cite{K2, Lemmas 4.5 and 4.6}
for some special cases).
Since $1 \ge \|\pi(\vw)\|^{-u} $
whenever $u > \tfrac {v+1-{j}}{j}$ and $v $ is at least $ j-1$, we can conclude
that subgroups of rank greater than
$n-s$ have no impact on the \de\ of $\Cal L$. 

\medskip

It will be convenient to associate to $A$ the following quantities: for each $j = 1,\dots,n-s$,
define
$$
\omega_j(A) \df \sup\left\{ v\left| \aligned \exists\, 
\vw \in \Cal S_{n+1,{j}}\text{  with arbitrary large
   } \|\pi_{\bullet}(\vw)\|  \\
\text{   such that   }\|R_{\sssize A}
\vc(\vw)\|
< \|\pi_{\bullet}(\vw)\|^{-\frac {v+1-{j}}{j}}\ \ 
\endaligned\right.\right\}\,.\tag 5.2
$$
By Lemma 5.1, 
(4.12) holds if and only if $\omega_j(A) \le v$ for all $j = 1,\dots,n -s$.
Thus Corollary 4.3 can be rewritten as

 \proclaim{Corollary 5.2} For $\Cal L$  parametrized by  {\rm (0.6)},
 $\omega(\Cal L)  = \max\big(\,n,\ \omega_j(A)_{j = 1,\dots,n -s}\big)$.
\endproclaim



We will refer to $\omega_j(A)$ as  the {\sl \de\ of $A$ of order\/} $j$. 
The reason for this terminology is the observation, essentially made 
 in \cite{K2}, that

 \proclaim{Lemma 5.3} $\omega_1(A) = \omega(A)$.
\endproclaim

\demo{Proof} Take $\vv = \pmatrix 
\vp
\\
\vq
\endpmatrix \in  \bz^{n+1}\nz = \Cal S_{n+1,1}$ in place of $\vw$, 
where  $\vp\in\bz^{s+1}$ and $\vq\in\bz^{n-s}$,  and observe 
that the only possible choice of $J$ in (5.1) and (4.3) is  $J = \varnothing$. 
It follows that $\vc(\vv) = \vv$,  
$R_{\sssize A}
\vc(\vw) = \vp + A\vq$ and $\pi_{\bullet}(\vv) = \vq$; hence the inequality 
in (5.2) coincides with (0.5).
 \qed\enddemo

In view of the lemma, the  estimate  (0.7$\ge$) stated in the Introduction
gets to be  a  special case of  Corollary 5.2 corresponding to $j=1$.
However it is worthwhile to note that  this inequality can be proved in an elementary 
way, and even more can be said:

 \proclaim{Lemma 5.4} Let $\Cal L$ be parametrized as in {\rm (0.6)\/}. 
Then for any 
$\,u < \omega(A)$ there exists an infinite subset $\Cal A$ of $\bz^{n+1}$ such that
  $$
|\vy\vq + p | < \|\vq\|^{-v}\text{  for all }\vy\in\Cal L\text{  and all  but finitely 
many }(p,\vq)\in \Cal A\,.\tag 5.3$$
\endproclaim

This not only proves the  lower bound  (0.7$\ge$), but also provides a way to
approximate {\it all\/} points of $\Cal L$ uniformly by a fixed sequence of integers. 

\demo{Proof} One knows that
for any $v < \omega(A)$ and  infinitely many
$\vq\in
\bz^{n-{s}}$ one can find $\vp = (p_0,p_1,\dots,p_s)\in\bz^{{s}+1}$ satisfying (0.5). Now take  any
$\vx\in \br^{s}$, denote $(p_1,\dots,p_s)$ by $\vp'$  and write
$$
\left|\,p_0 + (\vx,
\tilde\vx A) \left(\matrix \vp' \\ \vq
\endmatrix\right)\right|  = |p_0 + \vx\vp' + \tilde\vx A\vq| = |\tilde\vx( A\vq +
\vp)|\le 
\left\|\tilde\vx\right\|
\| A\vq + \vp\|\,.
$$ 
Slightly decreasing $v$ if needed, one gets 
$$
\left|\,p_0 + (\vx,
\tilde\vx  A) \left(\matrix \vp' \\ \vq
\endmatrix\right)\right| \le  \|\vq\|^{-v}
\tag 5.4
$$
for all but finitely many $(\vp,\vq)$ as above. Then it easily follows from (0.5) 
 that $\|\vp\|$ is bounded from above by $C\|\vq\|$,
where $C$ depends only on $A$. Thus, after possibly another slight change of $v$
and throwing away another finite subset, 
one can put $\left\|\pmatrix 
\vp'
\\
\vq
\endpmatrix\right\|^{-v}$ in the right hand side of (5.4).
\qed\enddemo

We now turn to the equality cases in   (0.7$\ge$), that is, to Theorem 0.4. The first observation
is that from Lemma 5.3
and Corollary 5.2 one can immediately deduce (0.7$=$) 
for $s = n-1$, that is, for column matrices $\Leftrightarrow$ codimension one   subspaces of $\br^n$. 
Since both $\omega(\Cal L)$ and $ \omega(A)$
are obviously infinite if columns of $A$ are linearly dependent over $\bq$, this proves 
Theorem 0.4 in case (b).
 
\medskip

Another case when  Theorem 0.4 holds for trivial reasons is 
$\Cal L = \{\va\}$,  a zero-dimensional subspace  represented by a $1\times n$ matrix $\va\in\br^n$. Indeed, it is a tautological statement that $\omega(\Cal L)$,
that is, the \de\ of the $\delta$-measure supported at $\va$, is equal to $\omega(\va)$.
On the other hand everything done in \S4 is easily applicable in the case $s = 0$ 
(and $\tilde\vx = 1$). Thus it follows from Corollary 5.2 that $\omega_j(\va)\le \omega(\va)$
for each $j$. (Exercise: prove it directly from the definition (5.2).)

To finish the proof of Theorem 0.4 it remains to treat the case  when $A$ is a matrix
with rationally proportional rows. For that it will be useful to get a better understanding of the 
`hidden symmetries' of  higher order exponents.  The three lemmas below serve this purpose.

 \proclaim{Lemma 5.5} For any $A\in M_{s+1,n-s}$ and all $\vw\in \Cal S_{n+1,{j}}$, 
$2\le j \le n-s$, one has 
$$
\max_{i = 0,\dots,s} \max\Sb J
\subset\{0,\dots,n\}\\ \# J = j-1
\endSb\big|\big\langle (\ve_i + \va_i)\wedge\ve_{\sssize J},\vw\big\rangle\big| \ll \|R_{\sssize A}
\vc(\vw)\| \tag 5.5a
$$
and
$$
 \|R_{\sssize A}
\vc(\vw)\| \ll \max_{i = 0,\dots,s} \max\Sb J
\subset\{i+1,\dots,n\}\\ \# J = j-1
\endSb\big|\big\langle (\ve_i + \va_i)\wedge\ve_{\sssize J},\vw\big\rangle\big|\,.
\tag 5.5b
$$
\endproclaim

Note that in both cases the reverse inequalities are obvious. The statement is also obvious 
for $j = 1$ (since, as was mentioned before, the only possible choice of $J$ is $J=\varnothing$).
Equation (5.5b)
allows one to slightly reduce  the set of pairs  $(i,J)$ involved in  computations of \de s.
Also, (5.5a) says that one can enlarge the set of pairs  in (5.1) so that
the formula, and hence the definition of $\omega_j(A)$, become symmetric 
under any\footnote{The invariance of 
 $\omega_j(A)$ under permutations
 not involving the top row  is an immediate consequence of (5.1); however the fact
that   the top row can also be permuted is nontrivial.} permutation of rows of $A$. 

\demo{Proof} A crucial observation is the following: for any $0\le k\le s$, using (5.1) 
and the linearity of $\vu\mapsto \langle (\ve_k + \va_k)\wedge\vu,\vw\big\rangle$,  
 one can write
$$
\big|\big\langle (\ve_k + \va_k)\wedge\vu,\vw\big\rangle\big|
\ll \|R_{\sssize A}
\vc(\vw)\|\tag 5.6a
$$
whenever $\vu\in\bigwedge^{j-1}(V_0)$ has norm $\ll \|A\|$. 
Likewise, assuming in addition that 
$\vu\in\bigwedge^{j-1}(\br\ve_{k+1}\oplus\dots\oplus\br\ve_{n})$, one has
$$
\big|\big\langle (\ve_k + \va_k)\wedge\vu,\vw\big\rangle\big|
\ll \max\Sb J
\subset\{k+1,\dots,n\}
\endSb\big|\big\langle (\ve_k + \va_k)\wedge\ve_{\sssize J},\vw\big\rangle\big|\,.\tag 5.6b
$$
To prove (5.5a),  take 
 $I\subset \{1,\dots,n\}$ of cardinality $j-2$, 
and  for any $0\le k\le s$ write
$$
\split
\big|\big\langle (\ve_k + \va_k)\wedge\ve_0\wedge\ve_{\sssize I},\vw\big\rangle\big| = 
&\ \big|\big\langle \ve_0\wedge(\ve_k + \va_k)\wedge\ve_{\sssize I},\vw\big\rangle\big|\\ \un{(5.6a) with 
$\vu = (\ve_k + \va_k)\wedge\ve_{\sssize I}$}\ll
 &\ \|R_{\sssize A}
\vc(\vw)\| + \big|\big\langle \va_0\wedge(\ve_k + \va_k)\wedge\ve_{\sssize I},\vw\big\rangle\big|
 \\  = 
 &\ \|R_{\sssize A}
\vc(\vw)\| + \big|\big\langle (\ve_k + \va_k)\wedge\va_0\wedge\ve_{\sssize I},\vw\big\rangle\big| 
\\ \un{(5.6a) with 
$\vu =  \va_0\wedge\ve_{\sssize I}$}\ll
 &\ \|R_{\sssize A}
\vc(\vw)\|\,.
\endsplit
$$
Similarly, if $k < i$ 
one has
$$
\split
\big|\big\langle (\ve_i + \va_i)\wedge\ve_k\wedge\ve_{\sssize I},\vw\big\rangle\big|= 
&\ \big|\big\langle \ve_k\wedge(\ve_i + \va_i)\wedge\ve_{\sssize I},\vw\big\rangle\big|\\ \un{(5.6b) with 
$\vu = (\ve_i + \va_i)\wedge\ve_{\sssize I}$}\ll
 \ \max\Sb J
\subset\{k+1,\dots,n\}
\endSb&\ \big|\big\langle (\ve_i + \va_i)\wedge\ve_{\sssize J},\vw\big\rangle\big| + 
\big|\big\langle \va_k\wedge(\ve_i + \va_i)\wedge\ve_{\sssize I},\vw\big\rangle\big|
 \\  \ll 
 \  \max\Sb J
\subset\{k+1,\dots,n\}
\endSb&\ \big|\big\langle (\ve_i + \va_i)\wedge\ve_{\sssize J},\vw\big\rangle\big| \,,
\endsplit
$$
and a repeated application of this trick allows one to reduce the estimation
of $\|R_{\sssize A}
\vc(\vw)\|$ to $J
\subset\{i+1,\dots,n\}$, proving (5.5b).   \qed\enddemo


The invariance of $\omega_j(A)$ with respect to permutations of rows of $A$
suggests that the same might hold for other row operations. This indeed happens to be the case:

 \proclaim{Lemma 5.6} Let $A' = BA$ for some $B\in \GL_{s+1}(\bq)$; in other words, $A'$ can 
be obtained from $A$  by a sequence of elementary row operations with rational coefficients
(that is, transposition of rows, addition of one row to another, and multiplication
of a row by a nonzero rational number). Then   $\omega_j(A') = \omega_j(A)$ for all $j$.
\endproclaim

\demo{Proof} Since row interchanges are taken care of in view of the previous lemma,
it remains to prove  $$\omega_j(A') \ge \omega_j(A)\tag 5.7$$  assuming
all the rows of $A'$ are the same as those of $A$ except for the top row,
and the latter is equal to:
{\rm (a)}  $\frac k\ell \va_0$, where $k$ and $\ell$ are nonzero integers, and
{\rm (b)} $ \va_0 + \va_1$.

For case (a), take $\vw\in \Cal S_{n+1,{j}}$,
write it in the form $\vw = \vw_0 + \ve_0\wedge \vw'$ where
both $\vw_0$ and $\vw'$ are in $\bigwedge(V_0)$, and put $$
\tilde \vw = \ell\vw_0 + k\ve_0\wedge \vw'\,.$$ 
It is easy to see  that  $
\tilde \vw$ also belongs to $\Cal S_{n+1,{j}}$, that is, represents a subgroup
of $\bz^{n+1}$: one can write 
$$
\vw = (a\ve_0 + b \vv_1)\wedge \vv_2\wedge\dots\wedge\vv_j
$$
for some integer vectors $\vv_1,\dots,\vv_j\in V_0$ and $a,b\in\bz$, and then take
 $$
\tilde \vw =(\ell a\ve_0 + kb \vv_1)\wedge \vv_2\wedge\dots\wedge\vv_j\,.$$
Now we claim that $\|R_{\sssize A'}
\vc(\tilde \vw)\|$ is not bigger than $\max\big(|k|,|\ell|\big)$ times $\|R_{\sssize A}
\vc(\vw)\|$. Indeed, one has, for  $J\subset \{1,\dots,n\}$,
$$
\split
\big\langle (\ve_0 + \tfrac k\ell\va_0)\wedge\ve_{\sssize J},\tilde\vw\big\rangle
 &
= \big\langle \ve_0 \wedge\ve_{\sssize J},k\ve_0\wedge \vw'\big\rangle  + 
\big\langle \tfrac k\ell\va_0 \wedge\ve_{\sssize J},\ell\vw_0 \big\rangle\\
&
= k\big( \big\langle \ve_{\sssize J}, \vw'\big\rangle  + 
\big\langle \va_0 \wedge\ve_{\sssize J},\vw_0 \big\rangle\big)
= k\,\big\langle (\ve_0 + \va_0)\wedge\ve_{\sssize J},\vw\big\rangle\,,
\endsplit
$$
and also, for $i\ge 1$,
$$
\big\langle (\ve_i + \va_i)\wedge\ve_{\sssize J},\tilde\vw\big\rangle
= \big\langle (\ve_i + \va_i)\wedge\ve_{\sssize J},l\vw_0\big\rangle
= \ell\,\big\langle (\ve_i + \va_i)\wedge\ve_{\sssize J},\vw_0\big\rangle= 
\ell\,\big\langle (\ve_i + \va_i)\wedge\ve_{\sssize J},\vw\big\rangle\,.
$$
This clearly implies (5.7). 
In case (b) the argument is similar. Namely, we take $\vw\in \Cal S_{n+1,{j}}$,
write it in the form $$\vw = \vw_0 + \ve_0\wedge \vw_0' + \ve_1\wedge \vw_1' +
\ve_0\wedge \ve_1\wedge \vw'\,,$$
 where
 $\vw_0,\vw',\vw_1',\vw_2' $ are all  in 
$\bigwedge(\br\ve_2\oplus\dots\oplus\br\ve_n)$, and put 
$$\tilde\vw = \vw + \ve_0\wedge  \vw_1'= \vw_0 + \ve_0\wedge (\vw_0' + \vw_1') + \ve_1\wedge \vw_1' +
\ve_0\wedge \ve_1\wedge \vw'\,.$$
Again, $
\tilde \vw$ can be easily shown to  represent a subgroup
of $\bz^{n+1}$:  write 
$$
\vw = (a\ve_0 + b \ve_1 + \vv_1)\wedge (c \ve_1 + \vv_2)\wedge\vv_3\wedge\dots\wedge\vv_j
$$
for some integer vectors $\vv_1,\dots,\vv_j\in \br\ve_2\oplus\dots\oplus\br\ve_n$ and $a,b,c\in\bz$, and then take
 $$
\tilde \vw =\big(a\ve_0 + b (\ve_0 + \ve_1) + \vv_1\big)\wedge \big(c (\ve_0 +  \ve_1) + \vv_2\big)\wedge\vv_3\wedge\dots\wedge\vv_j\,.$$
Now let us estimate $\|R_{\sssize A'}
\vc(\tilde \vw)\|$. It is clear that $\big\langle (\ve_i + \va_i)\wedge\ve_{\sssize J},\tilde\vw\big\rangle$ 
is the same as $\big\langle (\ve_i + \va_i)\wedge\ve_{\sssize J},\vw\big\rangle$
for    $i\ge 1$ and $J\subset \{1,\dots,n\}$. On the other hand,
$$
\split
\big\langle (\ve_0 + \va_0 + \va_1)\wedge\ve_{\sssize J},\tilde\vw\big\rangle
 &
= \langle \ve_0 \wedge\ve_{\sssize J},\vw  + \ve_0 \wedge\vw_1'\rangle  + 
\langle  \va_0  \wedge\ve_{\sssize J},\vw \rangle + 
\langle  \va_1  \wedge\ve_{\sssize J},\vw \rangle \\
&= \big\langle (\ve_0 + \va_0) \wedge\ve_{\sssize J},\vw\big\rangle + 
\langle \ve_{\sssize J},\vw_1'\rangle
+ 
\langle  \va_1  \wedge\ve_{\sssize J},\vw \rangle\\
&= \big\langle (\ve_0 + \va_0) \wedge\ve_{\sssize J},\vw\big\rangle + 
\big\langle (\ve_1 + \va_1) \wedge\ve_{\sssize J},\vw\big\rangle\,.
\endsplit
$$
Therefore $\|R_{\sssize A'}
\vc(\tilde \vw)\|\le 2\|R_{\sssize A}
\vc(\vw)\|$, so (5.7) holds as well. 
\qed
\enddemo

 \proclaim{Lemma 5.7} Suppose that $A$ has more than one row, and 
let $A'$ be the matrix obtained from $A$ by removing one of its
rows. Then  $\omega_j(A') \ge \omega_j(A)$ for all $j$.
If in addition the removed row  is a rational linear combination
of the remaining rows, then  $\omega_j(A') = \omega_j(A)$ for all $j$.
\endproclaim

\demo{Proof} 
Let $A'$ be obtained from $A$ by removing  its top
row (this can be assumed without loss of generality in view of the row 
interchange invariance).  It is clear that for any  $\vw\in \Cal S_{n+1,{j}}$  the inner product
  $\big\langle (\ve_i + \va_i)\wedge\ve_{\sssize J},\vw\big\rangle$ 
coincides with $\big\langle (\ve_i + \va_i)\wedge\ve_{\sssize J},\pi(\vw)\big\rangle$
whenever $i > 0$ and $ J
\subset\{1,\dots,n\}$. Thus it follows from (5.1) that $\big\|R_{\sssize A'}
\vc\big(\pi(\vw)\big)\big\|$ is not greater than $\|R_{\sssize A}
\vc(\vw)\|$, and obviously $\pi(\vw)\in \Cal S_{n,{j}}$ and $\pi_{\bullet}\big(\pi(\vw)\big) = \pi_{\bullet}(\vw)$.
Therefore whenever  $\vw\in \Cal S_{n+1,{j}}$ produces `a good approximation' to $A$
(meaning that $\|R_{\sssize A}
\vc(\vw)\|$ is smaller than $\|\pi_{\bullet}(\vw)\|$ to some negative power),
its projection $\pi(\vw)$ onto $\bigwedge^j(V_0)$  yields an equally good (or better) approximation
 to $A'$. This proves the first part of the lemma. 

For the second part, in view of Lemma 5.6 it is enough to assume that 
all the coefficients in the linear combination are zero, that is, $\va_0 = 0$.
Then one can reverse the above argument: whenever  
$\vw\in \bigwedge^j(V_0)$ produces `a good approximation' to $A'$,
it automatically yields an equally good  approximation
 to $A$, since one has $\big\langle (\ve_0 + \va_0)\wedge\ve_{\sssize J},\vw\big\rangle = 
\langle \ve_0 \wedge\ve_{\sssize J},\vw\rangle = 0$.
\qed
\enddemo

A combination of the above lemma with the observation made before Lemma 5.5
completes the proof of Theorem 0.4. In particular, we have shown that 
 subspaces $\Cal L$ of the form 
$$
\Cal L = \{(x_1,\dots,x_s,a_1,\dots,a_{n-s})\}
\quad\text{ 
or 
 }\quad\
\Cal L = \{(x_1,\dots,x_s,a_1x_i,\dots,a_{n-s}x_i)\}
$$
satisfy (0.7$=$), thus generalizing \cite{K2, Lemma 4.7} 
(whose method of proof, borrowed from \cite{BBKM},  did not shed any light 
on higher \de s of the corresponding matrices). Furthermore, in view of Lemma 5.4
one can conclude that whenever a subspace satisfying the assumptions of  Theorem 0.4
is not extremal (equivalently, $\omega(A) > n$), for any $u <  \omega(\Cal L) = \omega(A)$
one can find an infinite 
$\Cal A\subset\bz^{n+1}$ such that (5.3) holds; that is, 
there exists an infinite supply of  approximating vectors which can work uniformly for all points of $\Cal L$.

\heading{6. Generalizations and open questions}
\endheading

\subheading{6.1}
For matrices with no rational dependence between rows or columns 
the exponents of  orders higher than $1$ seem to be hard to understand. In particular the 
following question, a special case of which was asked in \cite{K2}, appears to be interesting:

\proclaim{Question} Does there exist a matrix $A\in M_{s+1,n-s}$ and  $2\le j \le n - s$
such that   $\omega_j(A)$ is
greater than both $n$ and $\omega(A)$?\endproclaim

An affirmative answer to this question would give a counterexample to (0.7$=$),
and, moreover, would provide an example of 
a proper nonextremal affine subspace $\Cal L$ of $\br^n$ such that for some
$\,v < \omega(\Cal L)$ it is impossible to find  an infinite subset 
$\Cal A$ of $\bz^{n+1}$ such that
$|\vy\vq + p | < \|\vq\|^{-v}$  for all $\vy\in\Cal L$  and all  but finitely 
many $(p,\vq)\in \Cal A$.
On the other hand, a useful consequence
of the validity of (0.7$=$) for all $\Cal L$ would be a possibility to compute the \hd\ 
of the set of subspaces of a given dimension and \de: indeed, then, in view of 
\cite{Do}, one
would have
$$
\dim\left(\left\{\Cal L \subset \br^n\left| \aligned \dim(\Cal L) = s\\ \omega(\Cal L) = v
\ \endaligned
\right.\right\}\right) 
= \cases ({s}+1)(n-{s}-1) + \frac{n+1}{v+1}\quad \text{if }v > n\,,\\ (s+1)(n-s)\quad\qquad
\qquad\text{otherwise. }
\endcases\tag 6.1
$$
By Theorem 0.4 this holds unconditionally for codimension one subspaces, that is,
$$
\dim\left(\left\{\Cal L \subset \br^n\left| \aligned \dim(\Cal L) = n-1\\ \omega(\Cal L) = v
\ \ \ \endaligned
\right.\right\}\right) 
= \cases  \frac{n+1}{v+1}\quad \text{if }v > n\,,\\ n
\qquad\text{otherwise. }
\endcases
$$
The author unfoundedly suspects  (6.1) to be true regardless of the answer to the above question;
 in other words,
even if higher order \de s can interfere with (0.7$=$), they conjecturally should not
be powerful enough to affect the computation of the \hd.

\subheading{6.2} The simplest matrices to look for possible counterexamples to  (0.7$=$) would be of size $2\times 2$,
corresponding  to lines in $\br^3$. To convince the reader
that the problem is far from trivial, let us work out an explicit formula for the second
order \de\ of $A = \pmatrix \va_0 \\ \va_1\endpmatrix= 
\pmatrix a_{02} & a_{03} \\ a_{12} & a_{13}\endpmatrix$. Write $\vw\in\Cal S_{4,2}$ in the form
$$
\vw = p\, \ve_0\wedge \ve_1 + \sum\Sb{i = 0,\,1}\\ {j = 2,\,3}\endSb w_{ij}\,\ve_i\wedge \ve_j 
+ q\,\ve_2\wedge \ve_3\,,
\tag 6.2
$$
and observe that for $i = 0,1$ one has
$$
\aligned
\big\langle (\ve_i + \va_i) \wedge\ve_{2},\vw\big\rangle = w_{i2} - a_{i3}q\,,\\
\big\langle (\ve_i + \va_i) \wedge\ve_{3},\vw\big\rangle = w_{i3} + a_{i2}q\,,
\endaligned\tag 6.3
$$
and also
$$
\aligned
\big\langle (\ve_0 + \va_0) \wedge\ve_{1},\vw\big\rangle &= p - a_{02} w_{12} - a_{03} w_{13}\\
&= p - \det (A) q - a_{02} (w_{12} - a_{13}q)  - a_{03}( w_{13}+ a_{12}q) \,.
\endaligned\tag 6.4
$$
Combining (6.3) and (6.4), one concludes  that the definition (5.2) of $\omega_2(A)$ reduces
to the following: $\omega_2(A)$ is the supremum of $v$ for which there exist
$\vw\in\Cal S_{4,2}$ of the form (6.2) with arbitrary large $|q|$ such that
$$
\max\big(|w_{02} - a_{03}q|,|w_{12} - a_{13}q|,  |w_{03} + a_{02}q|,  |w_{13} + a_{12}q|,
|p - \det (A) q|\big) <|q|^{-\frac{v-1}2}\,.
$$
Even when $\det (A) = 0$ (that is, rows/columns of $A$ are linearly dependent over $\br$ 
but {\it not\/} over $\bq$), the situation does not seem to be any less complicated.

\subheading{6.3} Theorem 2.2 can be 
used to 
treat  the so-called
{\it multiplicative\/} versions of the \di\ problems discussed in this paper.
Namely, define $$\Pi_{\sssize +}(\vq) 
\df \prod_{q_i \ne 0 } |q_i|\,,$$  denote by 
$\Cal W^\times_v$ the set of $\vy\in\br^n$ for which  there are infinitely many 
$\vq\in
\bz^n$ such that
$$
 |\vy\vq + p|   < \Pi_{\sssize +}(\vq)^{-v/n}  
$$
for some $p\in\bz$, and then define  {\sl multiplicative \de s\/}:

\roster
\item"$\bullet$" $\omega^\times(\vy)$ of $\vy$ by
$
\omega^\times(\vy) \df \sup\{v\mid \vy\in \Cal
W^\times_v\}\,,
$
\item"$\bullet$" $\omega^\times(\mu)$ of $\mu$ by
$
\omega^\times(\mu)\df 
\sup
\big\{\,v\bigm| \mu(\Cal
W^\times_v ) > 0\big\}\,.
$
\endroster
 
It is  easy to see that  $\omega^\times(\vy)$ is not less than $ \omega(\vy)$ for
all
$\vy$, and  $ \omega^\times(\vy) = n$ for $\lambda$-a.e.\ $\vy\in\br^n$,
that is, $\omega^\times(\lambda) = n$. 
To adapt the methods of the present paper to this set-up, one needs, following
 \cite{KM} and \cite{K2, \S 5}, to replace the one-parameter flow (3.2) by the action of the multi-parameter semigroup
$$
\text{\rm diag}(e^{t},
e^{-t_1},\dots,e^{-t_n})
\,,\quad\text{where }t_i \ge 0\text{ and } t = t_1 + \dots + t_n\,.
$$
This way it should be possible to prove Theorem 0.3 with $\omega$ replaced by $ \omega^\times$
(this was done in  \cite{K2} under the assumption that $\Cal L$ is {\sl strongly extremal\/},
that is, $\omega^\times(\Cal L) = n$),
and also derive formulas for  multiplicative \de s
of affine subspaces. 

\subheading{6.4}  With some abuse of notation, let us denote by $\sigma(\vy)$ the \de\ of 
the {\it column\/} matrix
given by the vector $\vy\in\br^n$, that is,
$$
\sigma(\vy) = \sup\big\{v \bigm| \exists\,\infty\text{-many }q\in\bz \text{ with }
\|q\vy + \vp\| < |q|^{-v}\text{ for some }\vp\in\bz^n\big\}\,,
$$
and, similarly to what was done with $\omega$, extend it to measures 
on $\br^n$.  It follows from Khintchine's Transference Principle
that  $\omega(\vy) = n$ if and only iff $\sigma(\vy)$ 
attains its smallest possible value, 
i.e.\ is equal to $1/n$. 
Therefore in all the problems related to the extremality of manifolds/measures, it makes no
difference whether to interpret points of $\br^n$ as row vectors (linear forms) or as column 
vectors 
(the latter set-up was employed in \cite{KLW, KW1--2} and is usually referred to 
as {\it simultaneous approximation\/}, hence our choice of notation $\sigma$).
The situation is however different when \de s are bounded away from their critical values.
Indeed, \cite{C, Ch.\ V, Theorem IV} estimates $\sigma(\vy)$
in terms of $\omega(\vy)$ as follows:
$$
\frac{\omega(\vy) - n +1}{n}\ge\sigma(\vy)\ge\frac{1}{
n-1 +  n/{\omega(\vy)}}\,,
$$
and the above inequalities are known to be sharp.
Thus in general it is not possible to extract any information concerning $\sigma(\mu)$ when $
\omega(\mu)$ is known. 

On the other hand, the methods of this paper can be 
adapted to computations
of the `simultaneous \de s' of manifolds and measures. Indeed, given a
column vector  $\vy\in\br^n$ one simply needs to work with the collection of vectors  of the form 
$\left(\matrix q\vy + \vp \\  q \endmatrix
\right)
$,
where $\vp\in\bz^n$ and 
$q \in \bz$, which is the same as $u_\vy\bz^{n+1}$, where one uses
$
u_\vy \df \left(\matrix
I_n & \vy  \\
0 & 1
\endmatrix \right)
$
instead of (3.1). Then, to study the simultaneous approximation  properties of
$\vy$, one uses the action by
$$
 g_t = \text{\rm diag}(e^{t/n},\dots,e^{t/n},
e^{-t})
$$
 which  uniformly  expands the first $n$ coordinates of vectors in $\br^{n+1}$ 
and contracts the last one. This way an application of Theorem 2.2  can yield Theorem 0.3 
with $\omega$ replaced by $\sigma$, as well as, after a multi-parameter modification,
by its multiplicative analogue.

\subheading{6.5} Studying \de s of matrices is a special case of a more general problem,
where one  replaces the right hand side of (0.5) by an
arbitrary function of
$\|\vq\|$.
Let us specialize to the case of row
vectors and use the following definition \cite{K2, \S 6.3}: for a nonincreasing function
$\psi:\bn\to (0,\infty)$, let  $\Cal W_\psi$ stand for the set of 
$\vy\in\br^n$ for which  there are infinitely
many
$\vq\in
\bz^n$ such that
$$
 \|\vy\vq + p\|   \le \psi(\|\vq\|) \quad \text{for some
}p\in\bz\,. 
$$
By Groshev's Theorem (see e.g.\  \cite{Sc2})  almost
no (resp., almost all)
$\vy\in\br^n$ belong to
$\Cal W_\psi$ if the series
$$
\sum_{k = 1}^\infty {k^{n-1}\psi(k)}\tag 6.5
$$
converges (resp., diverges). More generally, in \cite{BD} submanifolds $\Cal M$ of $\br^n$ are called 
{\sl of Groshev type for convergence\/} (resp., {\sl divergence\/}) if the convergence 
(resp.,  divergence) of (6.5) implies that  almost
no (resp., almost all) points of $\Cal M$ are in $\Cal W_\psi$. It is known
\cite{BKM, Be, BBKM} that nondegenerate submanifolds $\Cal M$ of $\br^n$ are of Groshev type for 
both convergence and divergence. 

The situation is much less understood when $\Cal M$ as above is replaced by a proper 
affine subspace $\Cal L$ of $\br^n$. Clearly $\omega(\Cal L)$ must be equal to $n$ in order 
for $\Cal L$ to 
be of Groshev type for convergence (since the choice $\psi(k) = k^{-v}$ for $v > n$ makes 
the series (6.1) converge) but the converse is not likely to be true.  On the other hand, 
in the following two cases it
 has been proved that $\Cal L$ parametrized as in (0.6) is  of Groshev type for convergence
under the assumption that $\omega(A)$ is strictly less than $n$: when  $\Cal L$ is a line passing 
through the origin \cite{BBDD}, and when it is of codimension one \cite{G1}. 
In the former case  $\Cal L$ was also proved to be of Groshev type for divergence.
Note that both cases fall into the framework of Theorem 0.4, that is, admit a simple
formula (0.7$=$) for the \de\ of $\Cal L$. In view of Corollary 5.2 it seems natural to make the 
following 

 \proclaim{Conjecture} Let $\Cal L$ be  parametrized by  {\rm (0.6)},
and suppose that $\omega_j(A)$ is strictly less than $n$ for every $j$. Then 
$\Cal L$ is of Groshev type for 
both convergence and divergence. 
\endproclaim

The author expects this conjecture, as well as the multiplicative analogue of 
its convergence case, to be provable by a combination of the methods 
of \cite{BKM, BBKM, G1,G3}
and the 
present paper.

\medskip

More generally, for any submanifold  of $\br^n$ it should be 
possible to state its own version of Groshev's Theorem, 
 with the convergence/divergence of (6.5) replaced by another `dividing line' condition. 
The following problems, posed in \cite{K2}, still  remain wide open:

\roster
\item"$\bullet$" is it true that the aforementioned `dividing line' condition of an affine
subspace
$\Cal L$ of
$\br^n$ is always inherited by manifolds nondegenerate in $\Cal L$? 
\item"$\bullet$" given a general affine subspace, say parametrized as in (0.6), 
find its `dividing line' condition, say in terms of the \di\ properties 
of the parametrizing matrix $A$; or, vice versa, 
describe the class of  subspaces with a given `dividing line'.
\endroster

\subheading{6.6} Finally, we remark that the generality of Theorem 2.2 
allows applications far beyond \da\ over $\br$. Namely, one can similarly consider
\di\ properties of measures on vector spaces over non-Archimedean local fields,
both of characteristic zero \cite{KT} and of positive characteristic \cite{G2}.
In fact, in the aforementioned two papers it was proved that manifolds nondegenerate
in the ambient spaces are extremal (and moreover strongly extremal). 
An application of Theorem 2.2 can extend these results to manifolds 
 nondegenerate in 
proper 
affine subspaces.


\Refs

\widestnumber\key{BBKM}



\ref\key Be \by V. Beresnevich  \paper A Groshev type theorem for convergence on
manifolds  \jour Acta Math. Hungar. \vol 94 \yr 2002\pages 99--130 \endref


\ref\key BBDD \by V. Beresnevich, V. Bernik, H. Dickinson, and M.\,M. Dodson
\paper On linear manifolds for which the Khinchin approximation theorem holds \jour
Vestsi Nats.
  Acad. Navuk Belarusi. Ser. Fiz.-Mat. Navuk  \yr 2000 \pages 14--17 
\lang Belorussian \endref


\ref\key BBKM \by V. Beresnevich, V. Bernik, D. Kleinbock, and G.\,A.
Margulis 
 \paper Metric Diophantine approximation: the Khintchine--Groshev theorem
for non-degenerate manifolds \jour Moscow Math. J. \yr
2002 \vol 2 \issue 2 \pages 203--225\endref


\ref\key BD \by V. Bernik and
M.\,M. Dodson \book Metric \da\ on
manifolds \publ Cambridge Univ. Press \publaddr Cambridge
\yr 1999 \endref



\ref\key BKM \by V. Bernik, D. Kleinbock, and G.\,A. Margulis \paper
Khintchine-type theorems  on
manifolds:  the convergence case for standard  and multiplicative
  versions \jour Internat. Math. Res. Notices \yr 2001   
\pages 453--486 \issue 9
\endref


\ref\key C \by J.\,W.\,S. Cassels \book An introduction to \di\ approximation \bookinfo
Cambridge Tracts in Math. \vol 45
\publ Cambridge Univ. Press \publaddr Cambridge
\yr 1957 \endref
 



\ref\key Da1 \by S.\,G. Dani \paper On invariant measures, minimal
sets, and a lemma of Margulis  \jour Invent. Math. \yr 1979
\pages 239--260 \issue 51
\endref


\ref\key Da2 \bysame \paper On orbits of unipotent flows on
homogeneous spaces \jour Ergod. Th. Dynam. Sys. \yr 1984 \pages
25--34 \issue 4
\endref



\ref\key Da3\bysame \paper Divergent trajectories of flows on
\hs s and Diophantine approximation\jour
J. Reine Angew. Math.\vol 359\pages 55--89\yr 1985\endref

\ref\key Da4 \bysame \paper On orbits of unipotent flows on
homogeneous spaces, II \jour Ergod. Th. Dynam. Sys. \yr 1986 \pages
167--182 \issue 6
\endref


\ref\key Do  \by M.\,M. Dodson \paper 
Hausdorff dimension, lower order and Khintchine's theorem in metric Diophantine
approximation
\jour  
J. Reine Angew. Math. \vol 432 \yr 1992 \pages 69--76\endref
 






\ref\key G1 \by A. Ghosh \paper A  Khintchine-type theorem for hyperplanes \jour J. London Math. 
Soc.  \vol 72 \issue 2  \yr 2005 \pages 293--304 \endref

\ref\key G2 \bysame \paper Metric \da\ over a local field of 
positive characteristic \jour J. Number Theory \vol 124 \yr 2007 \issue 2 \pages 454--469 \endref

\ref\key G3 \bysame \book Dynamics on homogeneous spaces and Diophantine approximation on
manifolds \bookinfo Ph.\,D. Thesis \publ Brandeis University \yr 2006 \endref



\ref\key {K1} \by D. Kleinbock \paper Some applications of
homogeneous dynamics to number theory \inbook in: Smooth Ergodic
Theory and Its Applications (Seattle, WA, 1999) \pages 639--660
\bookinfo Proc. Symp. Pure Math. \vol 68  \publ Amer. Math. Soc.
\publaddr Providence, RI \yr 2001 \endref

\ref\key K2 \bysame  \paper Extremal subspaces and their
submanifolds \jour Geom. Funct. Anal. \vol 13 \yr 2003 \issue 2 \pages
437--466 \endref

\ref\key K3 \bysame \paper Baker-Sprind\v zuk conjectures for complex analytic manifolds
\inbook in: Algebraic groups and Arithmetic  \publ TIFR, India \yr 2004 \pages 539-553
\endref






\ref\key KLW \by D. Kleinbock, E. Lindenstrauss,  and B. Weiss
\paper On fractal measures and \da
\jour Selecta Math. \vol 10 \issue 4 \pages  479--523 \yr 2004
         \endref

\ref\key KM \by D. Kleinbock and G.\,A. Margulis \paper Flows  on
homogeneous spaces and \da\ on manifolds\jour Ann. Math. \vol 148 \yr
1998 \pages 339--360 
 \endref




\ref\key KT \by D. Kleinbock and 
G. Tomanov \paper Flows on $S$-arithmetic homogeneous spaces and applications to metric
\da 
\jour Comm. Math. Helv. \vol 82 \pages 519--581 \yr 2007 \endref


\ref\key KW1 \by D. Kleinbock and B. Weiss \paper Badly approximable vectors on fractals
\jour 
Israel J.   Math. \vol 149  \yr 2005 \pages 137--170\endref

\ref\key KW2 \bysame \paper Friendly measures, homogeneous flows and singular vectors
\inbook in: 
Algebraic and Topological Dynamics,  Contemp. Math. \vol 385 \yr 2005 \publ AMS \publaddr
 Providence, RI \pages 281--292 \endref

\ref\key KW3 \bysame \paper Dirichlet's theorem on diophantine approximation and homogeneous flows
\jour J. Mod. Dyn. \vol 2 \yr 2008 \issue 1 \pages 43--62 \endref


\ref\key Mr1\by G.\,A. Margulis
\paper On the action of unipotent group in the space of lattices 
\inbook Proceedings of the Summer School on group representations  (Budapest 1971)\pages
365--370\publ Acad\'emiai Kiado
\publaddr Budapest \yr 1975\endref


\ref\key {Mr2}\bysame \paper Diophantine approximation, lattices and flows on homogeneous
spaces \inbook in: A panorama of number theory or the view from Baker's garden (Z\"
urich, 1999) \pages 280--310 \publ Cambridge Univ. Press \publaddr Cambridge\yr  2002 
\endref


\ref\key Mt\by P. Mattila \book Geometry of sets and measures
in Euclidean spaces. Fractals and rectifiability \bookinfo
Cambridge Studies in Advanced Mathematics, 44 \publ Cambridge
University Press \publaddr Cambridge \yr 1995 \endref

\ref\key PV \by A.  Pollington and S.  Velani \paper 
Metric Diophantine approximation and `absolutely friendly' measures
\jour Selecta Math.  \vol 11 \yr 2005 \issue 2 \pages 297--307\endref


\ref \key Rg \by M.\,S. Raghunathan \book Discrete subgroups of Lie
groups 
\publ Springer-Verlag \publaddr Berlin and New York \yr 1972 \endref%
 
\ref\key Rt1 \by  M. Ratner \paper Raghunathan's topological conjecture and
distributions of unipotent flows \jour  Duke Math. J. \vol 63 \yr 1991 \pages
235--280\endref

\ref\key Rt2 \bysame \paper Invariant measures and orbit closures for
unipotent actions on homogeneous spaces \jour Geom. Funct. Anal.\vol  4
\yr 1994\pages 236--257
         \endref






\ref\key Sc1\by W. Schmidt \paper Diophantine approximation and certain
sequences of 
lattices \jour Acta Arith. \vol 18 \yr 1971 \pages 195--178\endref

\ref\key Sc2 \bysame \book \di\ approximation \publ Springer-Verlag \publaddr Berlin and New York
\yr 1980 \endref%
 









\ref\key Sp\by V. Sprind\v zuk  \paper Achievements and problems in
Diophantine approximation theory \jour Russian Math. Surveys  \vol 35 \yr 1980 \pages 1--80 \endref


\ref\key SU \by  B. Stratmann and M. Urba\'nski \paper  Diophantine extremality of the
Patterson measure \jour Math. Proc. Cambridge Phil. Soc.  
\toappear 
\endref

\ref\key U \by  M. Urba\'nski \paper  Diophantine
approximation of self-conformal measures\jour J. Number Th. \vol 110 \yr 2005 \pages
219--235
\endref �
\endRefs
\enddocument